\documentclass[a4paper,11pt]{article}
\usepackage{latexsym}
\usepackage{amssymb}
\usepackage{amsmath}
\usepackage{textcomp}
\usepackage{mathrsfs}
\usepackage{eufrak}

\usepackage[T2A]{fontenc}
\usepackage[cp866]{inputenc}
\usepackage[english,russian]{babel}

\makeatother  \righthyphenmin=2 \oddsidemargin=0mm \topmargin=-1in \textwidth=16cm
\usepackage{graphicx}
\begin{document}
\begin{center}
\textbf{Solution of the unconditional extremum problem
for a liner-fractional integral functional
on a set of probability measures and its application
in the theory of optimal control of semi-Markov processes}\\

\textbf{Shnurkov P.V.}\footnote{National Research University
"Higher School of Economics"{},
Moscow Institute of Electronics and Mathematics,
Department of Applied Mathematics,
pshnurkov@hse.ru}
\end{center}

\parshape=1 1.45cm 13.5cm
In this paper, a new method for solving
the problem of optimal control of semi-Markov processes with finitely many states
is considered.
A new form of the assertion on an extremum of a liner-fractional integral functional
given on a set of probability measures is formulated and proved.
This form underlies the theorem of optimal control strategy for semi-Markov processes.
It is proved that the solution of the optimal control problem
for a semi-Markov process with finitely many states
is completely determined by the extremum properties
of the so-called test function of the liner-fractional integral functional
which is the control quality index.
At the same time, an explicit analytic representation was obtained for this test function
in terms of the initial probability characteristics of the semi-Markov model.

\textbf{Introduction.}

In this paper, we study two fundamental problems of the probability theory.
First, we consider the general problem of determining the global extremum
of liner-fractional integral functional given on a set of finite sets of probability measures.
We prove a new generalized and sharpened form of the assertion on the extremum of such a functional.
In what follows, this assertion underlies the theorem on the optimal control strategy
for a semi-Markov process with finitely many states.
We note that these results were first formulated without proof in~\cite{10}.

We begin with studying the unconditional extremum problem for a liner-fractional integral functional.
First, we make some remarks on the bibliography.

In the optimization theory, there is a field of studies, where the goal functional
of the extremum problem is the ratio of two linear functionals
and the constraints are of linear character.
This field of studies is called the liner-fractional programming.

In this optimization region, there is a vast scientific literature.
The contemporary theory of this field of scientific investigations
is described in the fundamental monograph~\cite{1},
where theoretical results about solutions of the corresponding extremum problems
are given together with numerical methods for determining such solutions.
A detailed bibliography of studies in the field of liner-fractional programming
can also be found there.
We also mention some significant works in the recent years,
where theoretical and computational problems in this field of scientific investigations
are described~\cite{2},~\cite{3},~\cite{4}.

A special section of liner-fractional programming is formed by extremum problems,
where the goal functional is the ratio of two integrals.
The integrand in these integrals is assumed to be known, and the integrals
are taken with respect to the probability measure belonging to a certain set of probability measures
defined on a give measurable space.
A solution of the problem is a probability measure
on which such a functional attains a global extremum.
The functionals of this type can be called \textit{liner-fractional integral functionals}.
Extremum problems for liner-fractional integral goal functionals
defined on the set of probability distributions in a finite-dimensional space
were studied by V.A.~Kashtanov~\cite{5},~\cite{6}.
His results in the theory of unconditional extremum for such functionals
are presented most completely in Chapter~10 of the collective monograph~\cite{6}.
The main result of these studies is that the unconditional extremum of such functionals
is attained on degenerate probability distributions with a single point of growth.
But this result was obtained in~\cite{6} under very strict conditions
the main of which was the assumption that the goal functional has an extremum,
i.e., the original problem has a solution.
As was already noted, in the present paper, we present a new solution of
the unconditional extremum problem for a liner-fractional integral functional,
which significantly generalizes an sharpens the results obtained in~\cite{6}.
The principal distinction between the result obtained in this paper from the preceding results
is that, in the main assertion, we present the conditions under which
an extremum of the liner-fractional integral functional exists
and is attained on a degenerate distribution concentrated at a single point.
In this case, the point of concentration of the whole probability measure
is the point of global extremum of a function for which the explicit analytic representation
was obtained.
\bigskip

\textbf{1. Solution of the extremum problem for a liner-fractional integral functional.}

First, we specify the notion of liner-fractional integral functional
and the statement of the main extremum problem.

Let $(U, \mathscr{B})$ be the main measurable space.
We immediately note that, in the second part of this study,
this space is interpreted as the set of admissible solutions or controls by a certain stochastic process.
In the first part, such an interpretation is not necessary, and the pair $(U,\mathscr{B})$
can be understood as an abstract mathematical object.
We assume that the $\sigma$-algebra $\mathscr{B}$ contains all sets consisting of a single point,
i.e., any element $u \in U$ treated as a single set $u = \{u\}$
belongs to this $\sigma$-algebra, i.e., the set $\{u\} \in  \mathscr{B}$ for any $u \in U$.

We consider a set $\Gamma$ consisting of probability measures defined on the
$\sigma$-algebra $\mathscr{B}$. The elements of the set~$\Gamma$,
which we denote by the symbol~$\Psi$, $\Psi \in \Gamma$,
are real functions satisfying the Kolmogorov axioms.
The construction and the properties of probability measures
are completely described in A.N.~Shiryaev's book~(\cite{7}, Chapter~2).
The theory of probability measures is also successively and profoundly presented
in the monograph of P.~Hennequin and A.~Tortrat~\cite{8}.
The conditions imposed on the set $\Gamma$ will be formulated and
analyzed below.

Now we introduce a special notion of degenerate probability measure
which is often used in what follows.
\medskip

\textbf{Definition~1.}
A probability measure $\Psi^*$ defined on the $\sigma$-algebra $\mathscr{B}$
is said to be \textit{degenerate} if there exists an element $u^* \in U$
satisfying the conditions $\Psi(\{u^*\})=1$ and $\Psi^*(U \setminus \{u^*\})=0$,
where $\{u^*\}=u^*$ is the set consisting of a single point $u^* \in U$.
The corresponding point $u^* \in U$ will be called a \textit{concentration point}
of the degenerate probability measure $\Psi^*$.
Thus, each degenerate probability measure $\Psi^*$ is determined by its concentration point~$u^*$.
In what follows, we use the notation $\Psi_{u^*}^{*}$ bearing in mind that
the degenerate probability measure $\Psi^*$
is concentrated at the point~$u^*$.
We also note that the degenerate probability measure $\Psi_{u^*}^{*}$
corresponds to a deterministic quantity which takes the fixed value $u=u^*$
with probability~$1$.
\medskip

By $\Gamma^*$ we denote the set of all degenerate probability measures
defined on the measurable space ($U,\mathscr{B}$).
The set~$\Gamma^*$ is in the one-to-one correspondence with the set of concentration points
of probability measures, i.e., with the set~$U$.

Consider functions $A(u)$ and $B(u)$ defined on the set $U$ and
taking real values $A:U\to R$, $B:U\to R$.
We introduce mappings $I_1: \Gamma \to R$ and $I_2 : \Gamma \to R$ defined as
\begin{equation*}
I_1(\Psi) = \int\limits_{U} A(u) d\Psi(u), \quad I_2(\Psi) = \int\limits_{U}B(u)d\Psi(u),
\eqno(1)
\end{equation*}
where $\Psi\in\Gamma$ is a probability measure.
The integrals in relation~(1) are general Lebesgue integrals
on the measurable space $(U,\mathscr{B})$ with respect to the probability measure $\Psi$.
In this case, the integrands $A(u)$, $B(u)$ can be interpreted as random variables
defined on the corresponding probability space $(U,\mathscr{B},\Psi)$,
and the Lebesgue integrals, as mathematical expectations of these random variables.
In the context of the control problem for semi-Markov processes, which is considered in this  paper,
the functions $A(u)$ and $B(u)$ have the meaning of conditional mathematical expectation
of some characteristics of this process defined under the condition
that the solution $u\in U$ is accepted.
The meaning of these functions is described in more detail in Section~5 below.

If the role of the measurable space $(U,\mathscr{B})$ is played
by the one-dimensional Euclidean space~$R$ with standard $\sigma$-algebra of Borel sets $\mathscr{B}$,
then the above-cited integrals can be understood as Lebesgue--Stieltjes integrals
over the probability measure defined by the distribution function~$\Psi$.
The theory of such integrals is given in (\cite{7}, Chapter~2).
It is assumed that these integrals in relations~(1) exist for any probability measure~$\Psi$
in the set~$\Gamma$.

Now we can introduce the notion of liner-fractional integral functional.
Assume that the linear functionals $I_1 : \Gamma \to R$ and $I_2 : \Gamma \to R$
defined by~(1) are given. In addition, we assume that
the functional $I_2(\Psi)$ does not vanish on the considered set of probability measures~$\Gamma$,
i.e., $I_2(\Psi) \neq 0, \Psi \in \Gamma$.
Let us consider the mapping $I = I(\Psi): \Gamma \to R$ defined by the relation
$$
I(\Psi) = \dfrac{I_1(\Psi)}{I_2(\Psi)} = \dfrac{\int\limits_{U}A(u)d\Psi(u)}{\int\limits_{U}B(u)d\Psi(u)}
\eqno(2)
$$
\medskip

\textbf{Definition 2.}
The mapping $I(\Psi): \Gamma \to R$ defined by relation~(2) will be called
a \textit{liner-fractional integral functional on the set of probability measures}~$\Gamma$.
\medskip

\textbf{Definition 3.} The function $C(u) = \dfrac{A(u)}{B(u)}$, $u \in U$,
will be called the \textit{test function} of the liner-fractional functional $I(\Psi)$
defined by relation~(2).

The unconditional extremum problem for a liner-fractional integral functional of the form~(2)
will be called the \textit{extremum problem}
$$
I(\Psi)\rightarrow \mbox{extr},
\eqno(3)
$$
$$
\Psi\in\Gamma.
$$
\medskip

To somewhat generalize the problem under study and to harmonize the further\break consideration,
we shall cosier the multidimensional analog of extremum problem~(3).
For this, we introduce necessary mathematical constructions.

We introduce a finite set of measurable spaces
$(U_1,\mathscr{B}_1),(U_2,\mathscr{B}_2),\dots,(U_N,\mathscr{B}_N)$,
where $U_i$ is a set of possible admissible controls and
$\mathscr{B}_i$ is the $\sigma$-algebra of subsets of the set~$U_i$
containing all single-point subsets
of the set~$U_i$, i.e., if $u_i \in U_i$, then $\{u_i\} \in \mathscr{B}_i$, $i=1,2,\dots,N$.
Let $\Gamma_i$ be a set of probability measures
$\Psi_i \in \Gamma_i$
defined on the $\sigma$-algebra $\mathscr{B}_i$,
and let $\Gamma_i^{*}$ be the set of all degenerate measures defined
on the $\sigma$-algebra $\mathscr{B}_i$,
and it is required that the condition $\Gamma_i^{*}\subset \Gamma_i$
be satisfied for all $i=1, 2,\dots, N$.
An arbitrary probability measure~$\Psi_i$ describes a random variable ranging in~$U_i$,
and a degenerate measure $\Psi_i^*$ concentrated at the point~$u_i^*$
corresponds to the deterministic quantity $u_i^*\in U_i$.
We assume that the corresponding structures are defined on all measurable decision spaces
$(U_1,\mathscr{B}_1), (U_2,\mathscr{B}_2), \dots, (U_N,\mathscr{B}_N)$.

Following the theory of products of measurable spaces presented in P.~Halmos's book
(\cite{9}, Chapter~VII), we consider the products of sets $U=U_1\times U_2\times \dots \times U_N$
and of the corresponding $\sigma$-algebras
$\mathscr{B}=\mathscr{B}_1\times\mathscr{B}_2\times\dots\times\mathscr{B}_N$.
On the measurable space $(U,\mathscr{B})$, one can define the probability measure
$\Psi=\Psi_1\times \Psi_2\times \dots \times\Psi_N$ which is called the product of measures
$\Psi_1,\Psi_2,\dots,\Psi_N$.
In what follows, for the measure $\Psi$ we use the multidimensional (vector) notation
$\Psi=(\Psi_1, \Psi_2,\dots,\Psi_N)$ treating it as an object defined by its components
$\Psi_1, \Psi_2,\dots,\Psi_N$.

Consider the set of probability measures defined on the measurable space $(U,\mathscr{B})$
in the sense cited above. We denote this set by the symbol~$\Gamma$,
and its elements, by the symbol $\Psi\in \Gamma$.
With regard to the above assumptions, the set $\Gamma$ can be interpreted as
the set of finite collections of probability measures:
$\Gamma=\{\Psi=\left(\Psi_1,\Psi_2,\dots,\Psi_N\right)$, $\Psi_i \in \Gamma_i$, $i=1,2,\dots, N\}$.
Necessary requirements to the set~$\Gamma$ and the corresponding sets $\Gamma_1,\Gamma_2,\dots,\Gamma_N$
will be given below.
By $\Gamma^*=\{\left(\Psi_1^*, \Psi_2^*,\dots,\Psi_N^*\right)$, $\Psi_i^*\in\Gamma_i^*$,
$i=1, 2,\dots, N\}$, we denote the set of all possible collections of degenerate probability measures
or, which is the same, the set of all degenerate probability measures defined
on the measurable space $(U,\mathscr{B})$.

Now we introduce real measurable functions
$A(u_1, u_2,\dots, u_N): U=U_1\times U_2\times\dots\times U_N\rightarrow R$ and
$B(u_1, u_2,\dots, u_N): U=U_1\times U_2\times \dots\times U_N\rightarrow R$
and consider the mapping $I(\Psi)=I(\Psi_1, \Psi_2,\dots,\Psi_N): \Gamma \rightarrow R$
defined by the relation
\begin{equation*}
\begin{aligned}
&I(\Psi)=I(\Psi_1,\Psi_2,\dots,\Psi_N)=
\dfrac{\int\limits_{U_1}\cdots\int\limits_{U_N}A(u_1,u_2,\dots,u_N)d\Psi_1(u_1) d\Psi_2(u_2)\dots
d\Psi_N(u_N)}
{\int\limits_{U_1}\cdots\int\limits_{U_N}B(u_1,u_2,\dots,u_N)d\Psi_1(u_1) d\Psi_2(u_2)\dots
d\Psi_N(u_N)},
\end{aligned}
\eqno(4)
\end{equation*}
In what follows, the mapping $I(\Psi): U\rightarrow R$
will be called a \textit{liner-fractional (multidimensional) integral functional}.
The function $C(u_1, u_2,\dots, u_N)=\dfrac{A(u_1,u_2,\dots,u_N)}{B(u_1,u_2,\dots,u_N)}$
will be called the \textit{test function of liner-fractional integral functional}~(4).

We shall investigate the unconditional extremum problem for the above-introduced
liner-fractional integral functional
$$
I(\Psi)=I(\Psi_1, \Psi_2,\dots,\Psi_N)\rightarrow extr,
$$
$$
\Psi=(\Psi_1, \Psi_2,\dots,\Psi_N)\in\Gamma.
\eqno(5)
$$
For the extremum problem~(5) to be well posed, we impose some conditions on the mathematical objects
determining this problem. Since for each problem of this form,
the integrands in the numerator and denominator $A(u_1, u_2,\dots, u_N)$, $B(u_1, u_2,\dots, u_N)$
are assumed to be prescribed, such conditions actually determine the character of the set of collections
of probability measures $\Gamma=\{(\Psi_1, \Psi_2,\dots,\Psi_N)$, $\Psi_i\in \Gamma_i$, $i=1,2,\dots, N\}$,
on which the original problem~(5) is posed. We formulate these conditions in the following form.

\begin{enumerate}
\item[1.] The integral expressions
$$
I_1(\Psi)=I_1(\Psi_1,\Psi_2,\dots,\Psi_N)
=\int\limits_{U_1}\cdots\int\limits_{U_N}A(u_1,u_2,\dots,u_N)d\Psi_1(u_1) d\Psi_2(u_2)\dots
d\Psi_N(u_N)
$$
$$
I_2(\Psi)=I_2(\Psi_1,\Psi_2,\dots,\Psi_N)
=\int\limits_{U_1}\cdots\int\limits_{U_N}B(u_1,u_2,\dots,u_N)d\Psi_1(u_1) d\Psi_2(u_2)\dots
d\Psi_N(u_N)
$$
are defined for all considered for all sets of probability measures $\Psi=(\Psi_1,\Psi_2,\dots,\Psi_N)\in\Gamma$.

\item[2.] The functional $I_2(\Psi)=I_2(\Psi_1, \Psi_2,\dots,\Psi_N)\neq 0$
for all $\Psi=(\Psi_1,\Psi_2,\dots,\Psi_N)\in \Gamma$.

\item[3.] The set $\Gamma$ contains the set of all degenerate probability measures:
$\Gamma^* \subset \Gamma$.
\end{enumerate}

These conditions will be called basic preliminary conditions
which are briefly denoted by BPC.

We make an important remark about these preliminary conditions.
\medskip

\textbf{Remark~1.} It follows from condition~2 that the function $B(u_1,u_2,\dots,u_N)$
cannot take values of opposite signs.
With regard to condition~3 we see that this function must be of constant sign on the whole set~$U$.
On the other hand, if the function
$B(u_1, u_2,\dots, u_N),~(u_1, u_2,\dots, u_N)\in U$ is strictly of constant sign,
then condition~2 is also satisfied. This fact directly follows from the properties of integral.
\medskip

Now we formulate and prove the main assertion about the extremum
of the liner-fractional integral functional.
\medskip

\textbf{Theorem~1.} {\it Consider the extremum problem~(5)
for the liner-fractional integral\break functional~(4).
Assume that conditions~1 and~3 in the system of preliminary conditions~$BPC$
are satisfied. Assume also that the function $B(u_1,u_2,\dots,u_N)$
in the definition of functional~(4) is strictly of constant sign
(strictly positive or strictly negative) for all value of the arguments
$(u_1,u_2,\dots,u_N)\in U$. Then the following assertions are true.
\begin{enumerate}
\item[(1)]
If the test function $C(u_1,u_2,\dots,u_N)=\dfrac{A(u_1,u_2,\dots,u_N)}{B(u_1,u_2,\dots,u_N)}$
is bounded above or below and attains its global extremum on the set
$U=U_1\times U_2\times\dots\times U_N$ (maximum or minimum),
then the corresponding extremum of the liner-fractional integral functional exists
and is attained on the set of degenerate probability measures
$\Psi^{*(+)}=\left(\Psi_1^{*(+)},\Psi_2^{*(+)},\dots,\Psi_N^{*(+)}\right)\in\Gamma^*$
or on the set of degenerate probability measures
$\Psi^{*(-)}=\left(\Psi_1^{*(-)},\Psi_2^{*(-)},\dots,\Psi_N^{*(-)}\right)$
depending on the form of extremum (maximum or minimum).
In this case, the degenerate probability measures
$\Psi_1^{*(+)},\Psi_2^{*(+)},\dots,\Psi_N^{*(+)}$ are concentrated at the points
$u_1^{(+)},u_2^{(+)},\dots,u_N^{(+)}$,\break respectively, if
$u^{(+)}=\left(u_1^{(+)},u_2^{(+)},\dots,u_N^{(+)}\right)\in U$
is the point of global maximum of the function $C(u_1,u_2,\dots,u_N)$, and then
$$
\max\limits_{\Psi \in \Gamma} I(\Psi)=\max\limits_{\Psi_i \in \Gamma_i, i=\overline{1,N}} I(\Psi_1,\Psi_2,\dots,\Psi_N)=
\max\limits_{\Psi_i^* \in \Gamma_i^*, i=\overline{1,N}} I(\Psi_1^*,\Psi_2^*,\dots,\Psi_N^*)
$$
$$
=I\left(\Psi_1^{*(+)},\Psi_2^{*(+)},\dots,\Psi_N^{*(+)}\right)
=\max\limits_{(u_1,u_2,\dots,u_N)\in U}\dfrac{A(u_1,u_2,\dots,u_N)}{B(u_1,u_2,\dots,u_N)}
$$
$$
=\dfrac{A\left(u_1^{(+)}, u_2^{(+)},\dots,u_N^{(+)}\right)}{B\left(u_1^{(+)}, u_2^{(+)},\dots,u_N^{(+)}\right)}.
\eqno(6)
$$
Similarly, the degenerate probability measures $\Psi_1^{*(-)},\Psi_2^{*(-)},\dots,\Psi_N^{*(-)}$
are concentrated at the points $\left(u_1^{(-)}, u_2^{(-)},\dots,u_N^{(-)}\right)$, respectively, if
$u^{(-)}=\left(u_1^{(-)}, u_2^{(-)},\dots,u_N^{(-)}\right)\in U$ is the point of global minimum of the function
$C(u_1,u_2,\dots,u_N)$, and then
$$
\min\limits_{\Psi \in \Gamma} I(\Psi)=\min\limits_{\Psi_i \in \Gamma_i, i=\overline{1,N}} I(\Psi_1,\Psi_2,\dots,\Psi_N)=
\min\limits_{\Psi_i^* \in \Gamma_i^*, i=\overline{1,N}} I(\Psi_1^*,\Psi_2^*,\dots,\Psi_N^*)
$$
$$
=I\left(\Psi_1^{*(-)},\Psi_2^{*(-)},\dots,\Psi_N^{*(-)}\right)
=\min\limits_{(u_1,u_2,\dots,u_N)\in U}\dfrac{A(u_1,u_2,\dots,u_N)}{B(u_1,u_2,\dots,u_N)}
$$
$$
=\dfrac{A\left(u_1^{(-)}, u_2^{(-)},\dots,u_N^{(-)}\right)}{B\left(u_1^{(-)}, u_2^{(-)},\dots,u_N^{(-)}\right)}.
\eqno(7)
$$

\item[(2)]
If the test function $C(u_1,u_2,\dots,u_N)=\dfrac{A(u_1,u_2,\dots,u_N)}{B(u_1,u_2,\dots,u_N)}$
is bounded above or below but does not attain the global extremum on the set
$U=U_1\times U_2\times\dots\times U_N$, then for any $\varepsilon > 0$, one can choose
an $\varepsilon$-optimal set of degenerate probability measures
$\Psi^{*(+)}(\varepsilon)=\left(\Psi_1^{*(+)}(\varepsilon),\Psi_2^{*(+)}(\varepsilon),\dots,
\Psi_N^{*(+)}(\varepsilon)\right)$
or
$\Psi^{*(-)}(\varepsilon)=\left(\Psi_1^{*(-)}(\varepsilon),\Psi_2^{*(-)}(\varepsilon),\dots,
\Psi_N^{*(-)}(\varepsilon)\right)$\break
depending on the form of the extremum (maximum or minimum).
In this case, the degenerate probability measures
$\Psi_1^{*(+)}(\varepsilon), \Psi_2^{*(+)}(\varepsilon),\dots,\Psi_N^{*(+)}(\varepsilon)$
are concentrated at the points\break
$u_1^{(+)}(\varepsilon),u_2^{(+)}(\varepsilon),\dots,u_N^{(+)}(\varepsilon)$, respectively, where
$u^{(+)}(\varepsilon)=\left(u_1^{(+)}(\varepsilon),u_2^{(+)}(\varepsilon),\dots,u_N^{(+)}(\varepsilon)\right)\in
U$ is an arbitrary point such that
$$
\sup\limits_{(u_1,u_2,\dots,u_N) \in U}\dfrac{A(u_1,u_2,\dots,u_N)}{B(u_1,u_2,\dots,u_N)}-\varepsilon <
\dfrac{A\left(u_1^{(+)}(\varepsilon),u_2^{(+)}(\varepsilon),\dots,u_N^{(+)}(\varepsilon)\right)}
{B\left(u_1^{(+)}(\varepsilon),u_2^{(+)}(\varepsilon),\dots,u_N^{(+)}(\varepsilon)\right)}
$$
$$
<\sup\limits_{(u_1,u_2,\dots,u_N) \in U}\dfrac{A(u_1,u_2,\dots,u_N)}{B(u_1,u_2,\dots,u_N)}<\infty
\eqno(8)
$$
if the function $C(u_1,u_2,\dots,u_N)$ is bounded above and the extremum problem~(5) is the maximum problem.
Similarly, the probability measures
$\Psi_1^{*(-)}(\varepsilon),\Psi_2^{*(-)}(\varepsilon),\dots,\Psi_N^{*(-)}(\varepsilon)$
are concentrated at the points
$u_1^{(-)}(\varepsilon),u_2^{(-)}(\varepsilon),\dots,u_N^{(-)}(\varepsilon)$, respectively,
where\break
$u^{(-)}(\varepsilon)=\left(u_1^{(-)}(\varepsilon),u_2^{(-)}(\varepsilon),\dots,
u_N^{(-)}(\varepsilon)\right)\in
U$ is an arbitrary point such that
$$
-\infty<\inf\limits_{(u_1,u_2,\dots,u_N) \in U}\dfrac{A(u_1,u_2,\dots,u_N)}{B(u_1,u_2,\dots,u_N)} <
\dfrac{A\left(u_1^{(-)}(\varepsilon),u_2^{(-)}(\varepsilon),\dots,u_N^{(-)}(\varepsilon)\right)}
{B\left(u_1^{(-)}(\varepsilon),u_2^{(-)}(\varepsilon),\dots,u_N^{(-)}(\varepsilon)\right)}
$$
$$
<\inf\limits_{(u_1,u_2,\dots,u_N) \in U}\dfrac{A(u_1,u_2,\dots,u_N)}{B(u_1,u_2,\dots,u_N)}+\varepsilon
\eqno(9)
$$
if the function $C(u_1,u_2,\dots,u_N)$ is bounded below and the extremum problem~(5) is the minimum problem.
\item[(3)]
If the test function
$C(u_1,u_2,\dots,u_N)=\dfrac{A(u_1,u_2,\dots,u_N)}{B(u_1,u_2,\dots,u_N)}$
is not bounded above or below, then no optimal set of probability measures exists
in the sense of the corresponding extremum problem.
In this case, there exists a sequence of sets of degenerate probability measures
$\Psi^*(n)=\left(\Psi_1^{*(+)}(n),\Psi_2^{*(+)}(n),\dots,\Psi_N^{*(+)}(n)\right)$, $n=1,2,\dots$,
where the measures\break
$\Psi_1^{*(+)}(n),\Psi_2^{*(+)}(n),\dots,\Psi_N^{*(+)}(n)$ are concentrated at the points
$u_1^{(+)}(n),u_2^{(+)}(n),\dots,u_N^{(+)}(n)$, $n=1,2,\dots$, respectively,
for which the following relations are satisfied:
$$
I\left(\Psi^{*(+)}(n)\right)=I\left(\Psi_1^{*(+)}(n),\Psi_2^{*(+)}(n),\dots,\Psi_N^{*(+)}(n)\right)
$$
$$
=\dfrac{A\left(u_1^{(+)}(n),u_2^{(+)}(n),\dots,u_N^{(+)}(n)\right)}
{B\left(u_1^{(+)}(n),u_2^{(+)}(n),\dots,u_N^{(+)}(n)\right)}\to \infty~\mbox{ as } n\to\infty
\eqno(10)
$$
if the function $C(u_1,u_2,\dots,u_N)$ is not bounded above.
Similarly, there exist a sequence of sets of degenerate probability measures
$\Psi^{*(-)}(n)=\left(\Psi_1^{*(-)}(n),\Psi_2^{*(-)}(n),\dots,\Psi_N^{*(-)}(n)\right)$, $n=1,2,\dots$,
where the measures $\Psi_1^{*(-)}(n),\Psi_2^{*(-)}(n),\dots,\Psi_N^{*(-)}(n)$ are concentrated at the points
$u_1^{(-)}(n),u_2^{(-)}(n),\dots,u_N^{(-)}(n)$, $n=1,2,\dots$, respectively,
for which the following relations are satisfied:
$$
I\left(\Psi^{*(-)}(n)\right)=I\left(\Psi_1^{*(-)}(n),\Psi_2^{*(-)}(n),\dots,\Psi_N^{*(-)}(n)\right)
$$
$$
=\dfrac{A\left(u_1^{(-)}(n),u_2^{(-)}(n),\dots,u_N^{(-)}(n)\right)}
{B\left(u_1^{(-)}(n),u_2^{(-)}(n),\dots,u_N^{(-)}(n)\right)}\to -\infty~\mbox{ as } n\to\infty
\eqno(11)
$$
if the function $C(u_1,u_2,\dots,u_N)$ is not bounded below.
\end{enumerate}
In this case, the above-formulated assertions of each item in Theorem~1 can be satisfied
either separately for one of the two types of extremum or simultaneously for two types of extremum.}
\medskip

Before proving Theorem~1, we prove some auxiliary assertions.
\medskip

\textbf{Lemma~1.} {\it Consider a liner-fractional integral functional
$I(\Psi_1, \Psi_2,\dots,\Psi_N)$ of the\break form~(4) defined on a set of collections of probability measures
$\Psi=(\Psi_1,\Psi_2,\dots, \Psi_N) \in \Gamma$.
Assume that condition~1 in the $BPC$ set is satisfied on the set $\Gamma$
and the function $B(u_1, u_2,\dots, u_N)$ is strictly of constant sign for all
for all $(u_1, u_2,\dots, u_N) \in U$. Then the following assertions hold.

\begin{enumerate}
\item[(1)]
If the test function $C(u_1, u_2,\dots, u_N)=\dfrac{A(u_1,u_2,\dots,u_N)}{B(u_1,u_2,\dots,u_N)}$
is bounded above, i.e., if the condition
$$
C(u_1, u_2,\dots, u_N)=\dfrac{A(u_1, u_2,\dots, u_N)}{B(u_1, u_2,\dots, u_N)}\leq c_0^{(+)}<\infty , \quad
(u_1, u_2,\dots, u_N) \in U,
\eqno(12)
$$
is satisfied, then the inequality
$$
I(\Psi_1, \Psi_2,\dots, \Psi_N)\leq c_0^{(+)}
\eqno(13)
$$
holds for all $(\Psi_1, \Psi_2,\dots, \Psi_N) \in \Gamma$.
\item[(2)]
If the test function $C(u_1,u_2,\dots,u_N)=\dfrac{A(u_1,u_2,\dots, u_N)}{B(u_1,u_2,\dots,u_N)}$
is bounded below, i.e., if the condition
$$
C(u_1, u_2,\dots, u_N)=\dfrac{A(u_1, u_2,\dots, u_N)}{B(u_1, u_2,\dots, u_N)}\geq c_0^{(-)}>-\infty , \quad
(u_1, u_2,\dots, u_N) \in U,
\eqno(14)
$$
is satisfied, then the inequality
$$
I(\Psi_1, \Psi_2,\dots, \Psi_N)\geq c_0^{(-)}
\eqno(15)
$$
holds for all $(\Psi_1, \Psi_2,\dots, \Psi_N) \in \Gamma$.
\end{enumerate}
}

\textbf{Proof of Lemma~1.}
Let us prove the first assertion of the lemma. First, we assume that the function
$B(u_1,u_2,\dots,u_N)$ is strictly positive:
$$
B(u_1, u_2,\dots, u_N)>0,~~(u_1, u_2,\dots, u_N)\in U.
\eqno(16)
$$
We note that, in this case, by the property of the integral
(\cite{9}, Chapter~V),
$$
\int\limits_{U_1}\dots \int\limits_{U_N}B(u_1, u_2,\dots, u_N)
d\Psi_1(u_1)d\Psi_2(u_2)\dots d\Psi_N(u_N)>0
\eqno(17)
$$
for any fixed set $\Psi=(\Psi_1, \Psi_2,\dots, \Psi_N)\in \Gamma$.
From inequality~(12) with regard to~(16), we obtain
$$
A(u_1, u_2,\dots, u_N)\leq c_0^{(+)} B(u_1, u_2,\dots, u_N), (u_1, u_2,\dots, u_N)\in U.
\eqno (18)
$$
In turn, it follows from inequality (18) and the properties of the integral that
$$
\int\limits_{U_1}\dots \int\limits_{U_N}A(u_1, u_2,\dots, u_N) d\Psi_1(u_1)d\Psi_2(u_2)\dots d\Psi_N(u_N)
$$
$$
\leq c_0^{(+)}\int\limits_{U_1}\dots \int\limits_{U_N}B(u_1, u_2,\dots, u_N) d\Psi_1(u_1)d\Psi_2(u_2)\dots d\Psi_N(u_N)
\eqno(19)
$$
for any fixed set $\Psi=(\Psi_1,\Psi_2,\dots,\Psi_N)\in\Gamma$.
But then from~(19) with regard to~(17) we obtain
$$
I(\Psi_1, \Psi_2,\dots, \Psi_N)=\dfrac{\int\limits_{U_1}\dots
\int\limits_{U_N}A(u_1, u_2,\dots, u_N) d\Psi_1(u_1)d\Psi_2(u_2)\dots d\Psi_N(u_N)}
{\int\limits_{U_1}\dots \int\limits_{U_N}B(u_1, u_2,\dots, u_N) d\Psi_1(u_1)d\Psi_2(u_2)\dots
d\Psi_N(u_N)}\leq c_0^{(+)}
\eqno(20)
$$
for any fixed set $(\Psi_1, \Psi_2,\dots, \Psi_N)\in \Gamma$.

Now we assume that the function $B(u_1, u_2,\dots, u_N)$ is strictly negative:
$$
B(u_1, u_2,\dots, u_N)<0,~~(u_1, u_2,\dots, u_N)\in U .
\eqno(21)
$$

Then
$$
\int\limits_{U_1}\dots \int\limits_{U_N}B(u_1, u_2,\dots, u_N) d\Psi_1(u_1)d\Psi_2(u_2)\dots d\Psi_N(u_N)<0
\eqno(22)
$$
for any fixed set $(\Psi_1, \Psi_2,\dots, \Psi_N)\in \Gamma$.

As previously, we start from inequality~(12). If conditions (21) and (22) are satisfied,
then the character of inequalities~(18) and~(19) becomes opposite but the character of inequality~(20)
remains unchanged. Thus, for any function $B(u_1, u_2,\dots, u_N)$ that is strictly of constant sign,
condition~(12) implies inequality~(20) which coincides with~(13).
The first assertion of Lemma~1 is proved. The second assertion can be proved similarly.
The proof of Lemma~1 is complete.
\medskip

\textbf{Lemma 2.} {\it Consider a liner-fractional integral functional
$I(\Psi_1, \Psi_2,\dots,\Psi_N)$ of the\break
form~(4) defined on a set of collections of probability measures
$\Psi=(\Psi_1,\Psi_2,\dots, \Psi_N)\in \Gamma$. Assume that condition~1 in the $BPC$ set
is satisfied on the set $\Gamma$ and the function $B(u_1, u_2,\dots, u_N)$ is strictly of constant sign
for all $(u_1, u_2,\dots, u_N)\in U$. Then the following assertions hold.
\begin{enumerate}
\item[(1)]
If the test function $C(u_1, u_2,\dots, u_N)=\dfrac{A(u_1,u_2,\dots,u_N)}{B(u_1,u_2,\dots,u_N)}$
is bounded above but does not attain its maximal value, then the inequality
$$
I(\Psi_1, \Psi_2,\dots, \Psi_N)< \sup\limits_{(u_1, u_2,\dots, u_N)\in U} C(u_1, u_2,\dots, u_N)<\infty
\eqno(23)
$$
holds for all $(\Psi_1, \Psi_2,\dots, \Psi_N)\in \Gamma$.
\item[(2)]
If the test function $C(u_1,u_2,\dots,u_N)=\dfrac{A(u_1,u_2,\dots,u_N)}{B(u_1,u_2,\dots, u_N)}$
is bounded below but does not attain its minimal value,
then the inequality
$$
I(\Psi_1, \Psi_2,\dots, \Psi_N)> \inf\limits_{(u_1, u_2,\dots, u_N)\in U} C(u_1, u_2,\dots, u_N)>-\infty
\eqno(24)
$$
holds for all $(\Psi_1, \Psi_2,\dots, \Psi_N)\in \Gamma$.
\end{enumerate}
}

\textbf{Proof of Lemma~2.}
Let us prove the first assertion of the lemma. Since the set of values
of the test function $C(u_1,u_2,\dots, u_N)$ is bounded above, it has a finite upper bound:
$$
\exists \sup\limits_{(u_1, u_2,\dots, u_N)\in U} C(u_1, u_2,\dots, u_N)<\infty
$$
(see~\cite{11}, Chapter~1, \S 3, item~3.4, Theorem~1).

By assumption, the function $C(u_1, u_2,\dots, u_N)$ does not attain its maximal value.
So we have the inequality
$$
C(u_1, u_2,\dots, u_N)=\dfrac{A(u_1, u_2,\dots, u_N)}{B(u_1, u_2,\dots, u_N)}
< \sup\limits_{(u_1, u_2,\dots, u_N)\in U} C(u_1, u_2,\dots, u_N)<\infty,
\eqno(25)
$$
$$
(u_1, u_2,\dots, u_N)\in U.
$$
Starting from the strict inequality~(25), we proceed as in the proof of
Lemma~1 with respect to inequality~(12).
As a result, we obtain the strict inequality~(23).

The second assertion in Lemma~2 can be proved similarly.
The proof of Lemma~2 is complete.
\medskip

\textbf{Proof of Theorem~1.} We begin with the proof of assertion~(1).
First, we assume that the test function
$C(u_1, u_2,\dots, u_N)=\dfrac{A(u_1, u_2,\dots, u_N)}{B(u_1, u_2,\dots, u_N)}$
is bounded above and attains the global maximum on the set $U$ at a certain point
$u^{(+)}=\left(u^{(+)}_1,u^{(+)}_2,\dots,u^{(+)}_N\right)\in U$. Namely,
$$
\max\limits_{(u_1, u_2,\dots, u_N)\in U} C(u_1, u_2,\dots, u_N) =
C\left(u^{(+)}_1,u^{(+)}_2,\dots,u^{(+)}_N\right)<\infty.
$$
Then we have the relation
$$
C(u_1, u_2,\dots, u_N)=\dfrac{A(u_1, u_2,\dots, u_N)}{B(u_1, u_2,\dots, u_N)} \leq
C\left(u^{(+)}_1,u^{(+)}_2,\dots,u^{(+)}_N\right)<\infty,
\eqno(26)
$$
$$
(u_1, u_2,\dots, u_N)\in U.
$$
The conditions of Lemma~1 are satisfied and we can use its assertions.
By the first of them, if inequality~(26) is satisfies, then
$$
I(\Psi_1, \Psi_2,\dots, \Psi_N)\leq C\left(u^{(+)}_1,u^{(+)}_2,\dots,u^{(+)}_N\right)<\infty
\eqno(27)
$$
for all sets of probability measures $\Psi=(\Psi_1, \Psi_2,\dots, \Psi_N)\in \Gamma$ in question.

Thus, the set of values of the liner-fractional integral functional $I(\Psi_1,\Psi_2,\dots,\Psi_N)$
is bounded above for all $\Psi=(\Psi_1,\Psi_2,\dots,\Psi_N)\in \Gamma$.
Then there exists a finite upper bound of this set and the inequality
$$
\sup\limits_{(\Psi_1, \Psi_2,\dots, \Psi_N)\in \Gamma}
I(\Psi_1, \Psi_2,\dots, \Psi_N)\leq C\left(u^{(+)}_1,u^{(+)}_2,\dots,u^{(+)}_N\right)
\eqno(28)
$$
holds. We consider a specific set of probability measures
$\Psi^{*(+)}=\left(\Psi_1^{*(+)}, \Psi_2^{*(+)},\dots,\Psi_N^{*(+)}\right)$,
where each probability measure $\Psi_i^{*(+)}$ is degenerate and concentrated
at a point $u_i^{(+)}$, $i=\overline{1, N}$. By the property of the integral, we have
$$
I\left(\Psi_1^{*(+)}, \Psi_2^{*(+)},\dots, \Psi_N^{*(+)}\right)=C\left(u^{(+)}_1,u^{(+)}_2,\dots,u^{(+)}_N\right).
\eqno(29)
$$
From relations (28), (29), we obtain
$$
\sup\limits_{(\Psi_1, \Psi_2,\dots, \Psi_N)\in \Gamma} I(\Psi_1, \Psi_2,\dots, \Psi_N)\leq
I\left(\Psi_1^{*(+)}, \Psi_2^{*(+)},\dots, \Psi_N^{*(+)}\right).
\eqno(30)
$$
In addition, we note that the following belongingness relations are satisfied:
$$
\Psi^{*(+)}=\left(\Psi_1^{*(+)}, \Psi_2^{*(+)},\dots, \Psi_N^{*(+)}\right)
\in \Gamma^* \subset \Gamma ,
\eqno(31)
$$
From (31) and the properties of the upper bound we obtain
$$
\sup\limits_{\left(\Psi_1^{*}, \Psi_2^{*},\dots, \Psi_N^{*}\right)
\in \Gamma^*} I\left(\Psi_1^{*}, \Psi_2^{*},\dots, \Psi_N^{*}\right)
\leq \sup\limits_{\left(\Psi_1, \Psi_2,\dots, \Psi_N\right) \in \Gamma}
I\left(\Psi_1, \Psi_2,\dots, \Psi_N\right).
\eqno(32)
$$
We combine (29), (30) and, (32) to obtain
$$
\sup\limits_{\left(\Psi_1^{*}, \Psi_2^{*},\dots, \Psi_N^{*}\right) \in \Gamma^*}
I\left(\Psi_1^{*}, \Psi_2^{*},\dots, \Psi_N^{*}\right)\leq
\sup\limits_{\left(\Psi_1, \Psi_2,\dots, \Psi_N\right) \in \Gamma} I\left(\Psi_1, \Psi_2,\dots, \Psi_N\right)
$$
$$
\leq I\left(\Psi_1^{*(+)}, \Psi_2^{*(+)},\dots, \Psi_N^{*(+)}\right)
=\dfrac{A\left(u^{(+)}_1,u^{(+)}_2,\dots,u^{(+)}_N\right)}{B\left(u^{(+)}_1,u^{(+)}_2,\dots,u^{(+)}_N\right)}.
\eqno(33)
$$
From (33) with (31) taken into account, we derive that the maximum of the functional\break
$I(\Psi_1, \Psi_2,\dots,\Psi_N)$ on the set of collections of probability measures
$\Psi=(\Psi_1, \Psi_2,\dots, \Psi_N)\in \Gamma$ exists and is attained
on the specific set of degenerate probability measures\break
$\Psi^{*(+)}\left(\Psi_1^{*(+)}, \Psi_2^{*(+)},\dots, \Psi_N^{*(+)}\right)$ defined above.

Moreover, relation (33) directly implies relation~(6). Thus,
assertion~(1) is proved in the case where the test function $C(u_1, u_2,\dots, u_N)$
attains the global maximum. The corresponding assertion in the case where
the test function $C(u_1, u_2,\dots, u_N)$ attains the global minimum can be proved similarly.
In this case, we use the second assertion of Lemma~1.

Now we shall prove the second assertion of Theorem~1.
Assume that the test function
$C(u_1,u_2,\dots, u_N)=\dfrac{A(u_1, u_2,\dots, u_N)}{B(u_1, u_2,\dots, u_N)}$
is bounded above but does not attain the global maximum on the set
$U = U_1 \times U_2 \times \dots \times U_N$.
Then the set of values of the test function has a finite upper bound:
$$
C(u_1, u_2,\dots, u_N)=\dfrac{A(u_1, u_2,\dots, u_N)}{B(u_1, u_2,\dots, u_N)}
<\sup\limits_{(u_1, u_2,\dots, u_N)\in U} \dfrac{A(u_1, u_2,\dots, u_N)}{B(u_1, u_2,\dots, u_N)}<\infty,
\eqno(34)
$$
$$
(u_1, u_2,\dots, u_N)\in U.
$$
By the definition of upper bound, for any fixed $\varepsilon >0$, there exists a point\break
$(u_1^{(+)}(\varepsilon), u_2^{(+)}(\varepsilon),\dots, u_N^{(+)}(\varepsilon))\in U$
such that the double inequality~(8) is satisfied
(see~\cite{11}, Chapter~1, \S 3, item~3.4).
In other words, the value of the test function at this point lies in the left $\varepsilon$-neighborhood
of the upper bound.

We fix an arbitrary point
$(u_1^{(+)}(\varepsilon), u_2^{(+)}(\varepsilon),\dots,u_N^{(+)}(\varepsilon))\in U$
that has such a property and consider a specific set of degenerate probability measures
$\Psi^{*(+)}(\varepsilon)=\left(\Psi_1^{*(+)}(\varepsilon),
\Psi_2^{*(+)}(\varepsilon),\dots\right.$,\break $\left.\Psi_N^{*(+)}(\varepsilon)\right)$, where
the degenerate measure $\Psi_i^{*(+)}(\varepsilon)$
is concentrated at a point $u_i^{(+)}(\varepsilon)$,\break $i=1,2,\dots,N$.

By the property of the integral, we have
$$
I\left(\Psi_1^{*(+)}(\varepsilon), \Psi_2^{*(+)}(\varepsilon),\dots,
\Psi_N^{*(+)}(\varepsilon)\right)=C(u_1^{(+)}(\varepsilon), u_2^{(+)}(\varepsilon),\dots,
u_N^{(+)}(\varepsilon)).
\eqno(35)
$$
From relation~(35) with this property of the test function taken into account, we obtain
$$
\sup\limits_{(u_1, u_2,\dots, u_N)\in U} \dfrac{A(u_1, u_2,\dots, u_N)}{B(u_1, u_2,\dots, u_N)}-\varepsilon
<I\left(\Psi_1^{*(+)}(\varepsilon), \Psi_2^{*(+)}(\varepsilon),\dots, \Psi_N^{*(+)}(\varepsilon)\right)
$$
$$
<\sup\limits_{(u_1, u_2,\dots, u_N)\in U} \dfrac{A(u_1, u_2,\dots, u_N)}{B(u_1, u_2,\dots, u_N)}<\infty.
\eqno(36)
$$
We also note that, in the case under study, the conditions of Lemma~2 are satisfied.
We use the first assertion of this lemma, namely, relation~(23):
$$
I(\Psi_1, \Psi_2,\dots, \Psi_N)< \sup\limits_{(u_1, u_2,\dots, u_N)\in U}
\dfrac{A(u_1, u_2,\dots, u_N)}{B(u_1, u_2,\dots, u_N)}<\infty
\eqno(37)
$$
for all $(\Psi_1, \Psi_2,\dots, \Psi_N)\in\Gamma$.

It follows from relations (36) and (37) that the specific set of degenerate probability measures
$\Psi^{*(+)}(\varepsilon)=\left(\Psi_1^{*(+)}(\varepsilon), \Psi_2^{*(+)}(\varepsilon),\dots,
\Psi_N^{*(+)}(\varepsilon)\right)$ defined above is $\varepsilon$-optimal.

The second part of assertion~(2) of Theorem~1 which is related to the properties of the lower bound
can be proved similarly.

We shall prove the third assertion of Theorem~1. Assume that the set of values
of the test function $C(u_1,u_2,\dots, u_N)=\dfrac{A(u_1,u_2,\dots, u_N)}{B(u_1, u_2,\dots, u_N)}$
is not bounded above on the set $U=U_1\times U_2 \times \dots \times U_N$.
Then there exists a sequence of points
$\left(u_1^{(+)}(n), u_2^{(+)}(n),\dots,u_N^{(+)}(n)\right)\in U$, $n=1,2,\dots~$, such that
$$
C\left(u_1^{(+)}(n), u_2^{(+)}(n),\dots,u_N^{(+)}(n)\right)
=\dfrac{A\left(u_1^{(+)}(n), u_2^{(+)}(n),\dots,u_N^{(+)}(n)\right)}
{B\left(u_1^{(+)}(n), u_2^{(+)}(n),\dots,u_N^{(+)}(n)\right)}\longrightarrow \infty,
\eqno(38)
$$
$$
n\rightarrow \infty.
$$
We fix a certain sequence of points
$\left(u_1^{(+)}(n),u_2^{(+)}(n),\dots,u_N^{(+)}(n)\right)\in U$, $n=1,2,\dots$,
that have the above-cited property and consider the sequence
$\Psi^{*(+)}(n)=\left(\Psi_1^{*(+)}(n),\Psi_2^{*(+)}(n),\dots\right.$,\break
$\left.\Psi_N^{*(+)}(n)\right)$,
$n=1,2,\dots$, such that, for each $n$, the $\Psi_i^{*(+)}(n)$ is concentrated at the point
$u_i^{(+)}(n)$, $i=1, 2,\dots, N$, $n=1,2,\dots~$.
By the property of the integral, for any fixed $n=1,2,\dots$, we have
$$
I \left(\Psi^{*(+)}(n)\right)=I\left(\Psi_1^{*(+)}(n), \Psi_2^{*(+)}(n),\dots, \Psi_N^{*(+)}(n)\right)
$$
$$
=\dfrac{A\left(u_1^{(+)}(n), u_2^{(+)}(n),\dots,u_N^{(+)}(n)\right)}{B\left(u_1^{(+)}(n), u_2^{(+)}(n),\dots,u_N^{(+)}(n)\right)}.
\eqno(39)
$$
It follows from relations (38) and (39) that
$$
I\left(\Psi^{*(+)}(n)\right)=I\left(\Psi_1^{*(+)}(n), \Psi_2^{*(+)}(n),\dots, \Psi_N^{*(+)}(n)\right)\longrightarrow\infty,
\quad n \rightarrow\infty.
\eqno(40)
$$
Relation (40) means that the set of values of the liner-fractional integral functional\break
$I(\Psi_1, \Psi_2,\dots, \Psi_N)$ of the form~(4) is not bounded above on the set of collections
of degenerate probability measures
$\left(\Psi_1^{*(+)}(n), \Psi_2^{*(+)}(n),\dots,\Psi_N^{*(+)}(n)\right)\in\Gamma^*$,
and hence on a wider set of collections of probability measures
$(\Psi_1, \Psi_2,\dots, \Psi_N)\in\Gamma$.
In this case, the extremum problem~(5) in the form of the maximum problem
does not have a solution.
The corresponding assertion for the version where the set of values of the test function
$C(u_1,u_2,\dots,u_N)=\dfrac{A(u_1,u_2,\dots,u_N)}{B(u_1,u_2,\dots,u_N)}$
is not bounded below, can be proved similarly.
The third assertion of Theorem~1 is proved.
This completes the proof of Theorem~1.
\bigskip

\textbf{2. On the problem of optimal control of semi-Markov stochastic processes.}

The results obtained in the theory of extremum of a liner-fractional integral functional
can be applied to study the problem of optimal control of semi-Markov stochastic processes.
We assume that the mathematical statement of the problem of optimal control of semi-Markov processes
is the statement of unconditional extremum problem for a liner-fractional integral functional
(problem~(5)). Then the main assertions of Theorem~1 allow us to formulate
the corresponding assertions on the conditions for the existence
of an optimal control strategy for this process and the character of the optimal control.

Before we start to consider specific results in the theory of control
of semi-Markov stochastic processes
which are based on the theory of unconditional extremum of a liner-fractional integral functional,
we make several remarks about the bibliography. We stress that these remarks do not pretend to be
a complete survey of the theory of optimal control of Markov and semi-Markov processes.
The goal of these remarks is to determine the place of the results obtained in the present paper
between the known results of the theory of control of semi-Markov processes with finitely many states
and to establish their principal novelty and specific characteristics.

The theory of optimal control of Markov and semi-Markov stochastic processes
has been intensively developed since the 1960s.
Already in the first fundamental studies,
not only the problems of existence of optimal control strategies
but also methods for determining these strategies were considered.
Such problems with algorithmic content were solved by powerful methods of applied mathematics
which had been developed shortly before this, namely, linear programming and dynamic programming.
We first note R.~Howard's classical work~\cite{12}, where the method of dynamic programming
was used to solve problems of optimal control of continuous-time Markov processes.
Further, V.V.~Rykov~\cite{13} proved that the optimal strategy
in a similar model of Markov process control with overestimation taken into account
is also stationary.

An important role in the development of the theory of control of stochastic processes
was played by W.~Jewell's work~\cite{14},
where the semi-Markov models of control with and without regard to the overestimation
were considered.
This work was translated into Russian and was a starting point for many
subsequent investigations of domestic and foreign specialists.
In particular, B.~Fox showed~\cite{15} that the optimal control strategy
for semi-Markov processes without over\-estimation taken into account is a stationary strategy;
similar results were obtained by E.~Denardo for the version with overestimation~\cite{16}.

Among subsequent studies that have not only formally theoretical but also algorithmic content,
we mention those of R.~Howard~\cite{17}, B.~Fox~\cite{15}, and S.~Osaki and H.~Mine~\cite{18}.
In these works, the method of linear programming was used to determine
the optimal control strategies for semi-Markov processes.

The fundamental monograph by H.~Mine and S.~Osaki~\cite{19},
where the main results in the theory of optimal control of Markov and semi-Markov stochastic processes
were systematized, was published in 1970 and translated into Russian in 1977.
In fact, this book summarized the results of studies of stochastic control problems during ten years.
We note that the Markov and semi-Markov models of control with finite sets of states and
admissible solutions accepted in each state were considered in that monograph.
Some principal theoretical results were obtained; namely, it was proved that
the optimal control strategies for basic types of the considered models
with and without overestimation are deterministic and stationary.
Several procedures for determining the optimal control strategies
were developed and justified.
In particular, for the model of the semi-Markov process control without overestimation
when the set of states forms one ergodic class
and the control quality index is the stationary average specific income
(\cite{19}, Chapter~5, item~5.5),
the procedure for determining the optimal randomized strategy
was realized by the linear programming method.
We pay special attention to this result, because a similar model
of semi-Markov process control will be considered in the present paper.

A principal role in the development of the theory of stochastic control was played
was played by I.I.~Gikhman's and A.V.~Skorokhod's monograph~\cite{20},
where the foundations of the theory of optimal control of stochastic processes
with discrete and continuous time,
including the theory of control of processes described by stochastic differential equations,
were systematically presented for the first time.
The problems of control of Markov processes with discrete time
and jumpwise Markov processes with continuous time were considered separately.
The roles of sets of states and admissible controls were played by spaces of very general structure.
For wide classes of control quality functionals
(the so-called evolutionary functionals in Markov models with discrete time
and integral accumulation functionals in Markov models with continuous time),
theorems on the existence and forms of representation of optimal control strategies
were proved.
It was shown that, for homogeneous Markov models,
there exist optimal control strategies that are stationary and deterministic.
In other words,
such strategies are given by deterministic functions
whose argument is the system state at the time of acceptance of the solution
and which are independent of the time at which the solution is accepted.
As for the important question of the forms of representation of these functions,
they can be characterized as follows.
The functional equations complicated by the extremum condition satisfied by the above-cited functions
were determined.
In fact, these relations are the Bellman equations for the corresponding dynamic stochastic models.

We specially note that no control problems for semi-Markov processes
were considered in the monograph~\cite{20}.
But the further development of the general theory of control of such processes
followed the ideas outlined in that book.

In the further development of the theory of control of semi-Markov processes,
the models were complicated and the initial assumptions were generalized.
For example, in~\cite{21},~\cite{22}, the semi-Markov decision processes
were considered
under very general assumptions on the character of the state and decision spaces.
The control problems were studied for different types of goal indices
which generalize the above-cited stationary index of average specific profit.
It was proved in these papers that the optimal control strategy
with respect to each of the indices exists and is the same stationary deterministic strategy
given by a certain function defined on the set of states of the process.
It is only known about this function that it satisfies some integral equation
which is, in a sense, the Bellman equation for the corresponding control problem.

Among the works preceding this study we mention
V.A.~Kashtanov's paper published as chapter~13 of the collective monograph~\cite{6}.
In that paper, he considered the optimal control problem for
a semi-Markov process with finitely many states whose set of possible solutions
was the set of negative real numbers.
The model considered there is a model without overestimation, and
the control quality index is the stationary value of the average specific income
determined as in the classical works~\cite{14},~\cite{19}.
The randomized control in each state is determined by the probability distribution
introduced on the set of admissible solutions whose set is given by the control strategy.
The author stated that the stationary value of the average specific income
is in form a liner-fractional integral functional
of a set of probability distributions which form the control strategy.
In this case, it was previously shown (\cite{5};~\cite{6}, Chapter~10])
that the liner-fractional functional attains an extremum on degenerate distributions.
This naturally implies that the optimal control strategy
is deterministic and must be determined by an extremum point
of the function that is the ratio of the integrand in the numerator
to the integrand in the denominator of this liner-fractional functional.
But no explicit representations for the above-cited functions
were obtained in~\cite{6}.
Moreover, in the version of the extremum theorem for the liner-fractional integral functional,
which was given in Chapter~10 in the monograph~\cite{6},
it was required to verify the existence conditions for this extremum.
Such conditions were not given there.
Therefore, the results obtained in~\cite{6}
could not be used to prove the existence of a deterministic optimal control strategy
for a semi-Markov process and to obtain a rigorous justification of the method
for determining such a strategy.

Now we shall present the results of the optimal control theory
for a semi-Markov process obtained in the present study.
\bigskip

\textbf{3. Description of the mathematical model of control of a semi-Markov process.}

We consider the following model of control of a semi-Markov stochastic process
following the general scheme proposed in the classical works~\cite{14},~\cite{19}.
Let $\xi(t)$ be a semi-Markov stochastic process with finitely many states
$\{1,2,\dots, N\}$, $N<\infty$.
By $t_n$, $n=0,1,2,\dots$, $t_0=0$, we denote the random times of variation in the states of this process,
$\theta_n=t_{n+1}-t_n$, $n=0,1,2,\dots$, $\xi_n=\xi(t_n)=\xi(t_n+0)$, $n=0,1,2,\dots$
(we assume that the trajectory of the process $\xi(t)$ are right continuous).
The random sequence $\{\xi_n\}$ forms a Markov chain embedded in the semi-Markov process $\xi(t)$.
We specify the set of measurable spaces
$(U_1,\mathscr{B}_1),(U_2,\mathscr{B}_2),\dots,(U_N,\mathscr{B}_N)$,
where~$U_i$ is the set of admissible controls,
$\mathscr{B}_i$ is the $\sigma$-algebra of subsets of the set~$U_i$
including all single-point subsets of the set~$U_i$,
i.e., if $u_i \in U_i$, then $\{u_i\} \in \mathscr{B}_i$, $i=1,2,\dots,N$.
Assume that $\Gamma_i$ is a certain set of probability measures $\Psi_i \in \Gamma_i$
defined on the $\sigma$-algebra $\mathscr{B}_i$,
$\Gamma_i^{*}$ is the set of all degenerate probability measures defined on the $\sigma$-algebra~$\mathscr{B}_i$,
and the condition $\Gamma_i^{*}\subset \Gamma_i$ is satisfied for all $i=1,2,\dots,N$.

We assume that the controls of the semi-Markov stochastic process $\xi(t)$
are realized at times $t_n$, $n=0,1,2,\dots$, directly after variation in the state of the process.
If $\xi_n=\xi(t_n)=i \in X$, then the control is a random quantity~$u_n$
taking values in the set of admissible controls~$U_i$
and described by the probability measure (probability distribution) $\Psi_i \in \Gamma_i$.
We assume that, under the fixed condition $\xi_n=\xi(t_n)=i$, the control is determined
independently of the past behavior of the process $\xi(t)$
and the probability measure $\Psi_i$ describing the stochastic control~$u_n$
depends only on the state $i\in X$.
Then the choice of controls at the times of variation in the states $\{t_n, n=0,1,2,\dots \}$
is described by the set of probability measures (probability distributions)
$(\Psi_1, \Psi_2,\dots, \Psi_N)$, $\Psi_i \in \Gamma_i$, $i=1,2,\dots,N$.

Such a set will be called a \textit{control strategy} for the semi-Markov process $\xi(t)$.
In its properties, such a strategy is Markov, homogeneous, and randomized.

The control strategy $\left(\Psi_1^*, \Psi_2^*,\dots, \Psi_N^*\right)$, $\Psi_i^*\in\Gamma_i^*$,
$i=1,2,\dots,N$, defined by a set of degenerate probability measures
is said to be deterministic.

We shall use theoretical constructions related to a product of measurable spaces
considered in Section~1 above.
We let $\Gamma$ to denote the set of probability measures on
the $\sigma$-algebra $\mathscr{B}=\mathscr{B}_1\times \mathscr{B}_2\times \dots \times\mathscr{B}_N$
which are defined as a product of measures $\Psi_1,\Psi_2,\dots,\Psi_N$,
where $\Psi_i \in \Gamma_i$, $i=1,2,\dots,N$.
The set~$\Gamma$ can also be interpreted as a set of collections of
probability measures $\Psi=(\Psi_1,\Psi_2,\dots,\Psi_N)$, $\Psi_i\in\Gamma_i$, $i=1, 2,\dots, N$
or admissible control strategies for the semi-Markov process $\xi(t)$.
A more complete description of the set of admissible control strategies~$\Gamma$
and an analysis of conditions imposed on this set will be given in Section~5.

Further, by $\Gamma^*$ we denote the set of all possible collections of degenerate probability measures
$\Psi^*=\left(\Psi_1^*, \Psi_2^*,\dots, \Psi_N^*\right)$, $\Psi_i^*\in\Gamma_i^*$, $i=1, 2,\dots, N$
or all deterministic control strategies.

In what follows, we consider the control problem for the semi-Markov process $\xi(t)$,
where the control quality index has the form of a liner-fractional integral functional
which is defined on the set $\Gamma$ and can be represented in the form~(4).
Thus, the optimal control problem for a semi-Markov process $\xi(t)$
is formalized as the extremum problem~(5).

Further, we shall need different probability characteristics
of the semi-Markov decision process $\xi(t)$.
It is known from the general theory of semi-Markov processes
(\cite{25},~\cite{26}),
the main probability characteristic of such a process
is the so-called semi-Markov function. We define this function for a process with control
(\cite{19}, Chapter~5) as
$$
Q_{ij}(t,u)=P(\xi_{n+1}=j,\theta_n<t \mid \xi_n=i, u_n=u),~t\in [0,\infty),
u\in U_i;~i,j\in X=\{1,2,\dots,N\}.
\eqno(41)
$$
The use of semi-Markov functions allows one to obtain the transition probabilities
for the controlled Markov chain $\{\xi_n\}$ embedded in the semi-Markov process under study,
$$
p_{ij}(u)=P(\xi_{n+1}=j \mid \xi_n=i, u_n=u)= \lim\limits_{t\rightarrow\infty}Q_{ij}(t,u),~~
u\in U_i;~~ i,j\in X,
\eqno(42)
$$
and the distribution functions for the time intervals in which the semi-Markov process $\xi(t)$
stays in the corresponding states,
$$
H_{i}(t,u)=P(\theta_n<t \mid \xi_n=i, u_n=u)=\sum\limits_{j\in X}Q_{ij}(t,u),~~ t\in [0,\infty), u\in U_i;~~ i\in X.
\eqno(43)
$$

By
$$
T_{i}(u)=\textbf{E}\left[\theta_n \mid \xi_n=i, u_n=u\right]
=\int\limits_0^{\infty}\left[1-H_i(t,u)\right]dt,~~ u\in U_i;~~ i\in X,
\eqno(44)
$$
we denote the conditional mathematical expectations of the holding time
of the semi-Markov process $\xi(t)$ in each of the states.

The above-introduced characteristics (41)--(44) are defined under the conditions
that, at the time~$t_n$ of the state change, the process turns out to be in the state~$i$
and a solution $u\in U_i$ is accepted.
For a given control strategy $\Psi=\left(\Psi_1,\Psi_2,\dots,\Psi_N\right)$,
the corresponding probability characteristics can be written without conditions on the control.
Namely,
$$
Q_{ij}(t)=P(\xi_{n+1}=j,\theta_n<t \mid \xi_n=i)=\int\limits_{U_i}Q_{ij}(t,u) d\Psi_i(u),~~ t\in [0,\infty);~~ i,j\in X;
\eqno(45)
$$
$$
p_{ij}=P(\xi_{n+1}=j \mid \xi_n=i)=\int\limits_{U_i}p_{ij}(u) d\Psi_i(u),~~ i,j\in X;
\eqno(46)
$$
$$
T_{i}=\textbf{E}\left[\theta_n \mid \xi_n=i\right]=\int\limits_{U_i}T_{i}(u) d\Psi_i(u),~~ i\in X.
\eqno(47)
$$
In what follows, we assume that the probability characteristics
$p_{ij}(u)$, $u\in U_i$, $i,j\in X$, and $T_i(u)$, $u\in U_i$, $i\in X$,
defined by formulas (42) and (44) are given the semi-Markov model under study.
For the fixed control strategy $\Psi=(\Psi_1,\Psi_2,\dots,\Psi_N)$,
the corresponding probability characteristics $p_{ij}$, $T_i$, $i,j\in X$,
are defined by formulas~(46) and~(47) without conditions on the control.
\bigskip

\textbf{4. Stationary cost control quality index and its analytic structure.}

We introduce a certain additive cost functional related to the semi-Markov process~$\xi(t)$
under study. By its meaning, this functional is the random income or profit
accumulated during the time interval $[0,t]$.
Several definitions of such a functional are given in the fundamental works~\cite{14}, (\cite{19}, Chapter~5).
We let $\widetilde{v}(t)$, $t\geq 0$,
denote the value of this additive functional at time~$t$;
$\widetilde{v}_n=\widetilde{v}(t_n+0)$ is the corresponding value directly
after the next time $t_n$, $n=0,1,2,\dots$, of the state variation; and
$\widetilde{v}_0=v_0$ is a given value at time $t=0$.
We consider the quantity
$$
d_i(u)=\textbf{E}\left[\widetilde{v}_{n+1}-\widetilde{v}_n \mid \xi_n=i, u_n=u\right],~~u\in U_i,~~i\in X,
\eqno(48)
$$
which is the conditional mathematical expectation of
the increment of the additive cost functional in  the time interval between successive variations
in the states of the semi-Markov process~$\xi(t)$.
Then the corresponding mathematical expectation calculated without any condition on the solution
that is accepted at time~$t_n$ can be represented as
$$
d_i=\textbf{E}\left[\widetilde{v}_{n+1}-\widetilde{v}_n \mid \xi_n=i\right]=\int\limits_{U_i}d_i(u)d\Psi_i(u),~~i\in X.
\eqno(49)
$$

We assume that, for a given control strategy $\Psi=(\Psi_1,\Psi_2,\dots,\Psi_N)$,
the embedded Markov chain $\{\xi_n\}$ has precisely one class of recursive positive states
(by the terminology accepted in~\cite{19}, such a set of states is called an ergodic class).
As is known (\cite{7}, Chapter~VIII), this condition is necessary and sufficient for
the existence of a unique stationary distribution.
We denote this stationary distribution of the Markov chain $\{\xi_n\}$
by $\pi=(\pi_1, \pi_2,\dots, \pi_N)$.
We note that this distribution depends on the control strategy
$\Psi=(\Psi_1,\Psi_2,\dots,\Psi_N)$.
Under this condition, the following result holds,
which is called the ergodic theorem for the additive cost functional
$$
I=\lim\limits_{t\rightarrow\infty}\dfrac{\textbf{E}\widetilde{v}(t)}{t} =\dfrac{\sum\limits_{i=1}^N
d_i\pi_i}{\sum\limits_{i=1}^N T_i\pi_i}.
\eqno(50)
$$

Relation~(50) was proved in (\cite{19}, Chapter~5).
We note that similar results hold for much more general semi-Markov models (\cite{21},~\cite{22}).

From the standpoint of applications, the quantity defined by~(50) is the average specific profit
related to the system evolution in the stationary regime.

At the same time, according to the above remarks, the quantity~$I$
depends on the control strategy $\Psi=(\Psi_1,\Psi_2,\dots,\Psi_N)$.
Thus, the expression in the right-hand side of~(50) defines a functional on the set
of control strategies under study: $I(\Psi)=I(\Psi_1, \Psi_2,\dots, \Psi_N)$.

The applied meaning of the index $I(\Psi)$ determined its place in optimal control problems.
Starting from the first fundamental works in the theory of control
of semi-Markov processes~\cite{14}, \cite{19},
it played the role of control quality index or the goal functional.
Already in (\cite{19}, Chapter~5), it was mentioned that this functional has a liner-fractional form of dependence
on deterministic control strategies in the case of a finite set of admissible solutions
accepted in each state.
The next result related to the form of dependence of the functional $I(\Psi_1,\Psi_2,\dots,\Psi_N)$
on the control strategy $\Psi=(\Psi_1,\Psi_2,\dots,\Psi_N)$ was presented in the work~\cite{6}, Chapter~13,
characterized in Section~2.
Further, a stronger assertion was proved in~\cite{23} stating that the functional\break
$I(\Psi_1,\Psi_2,\dots,\Psi_N)$
analytically depends on the probability measures $\Psi_1,\Psi_2,\dots,\Psi_N$
that form the control strategy.
Now we completely formulate this assertion.

\textbf{Theorem 2.}
{\it The stationary cost index given by~(50) is a liner-fractional functional
of probability measures (probability distributions) $\Psi_1,\Psi_2,\dots,\Psi_N$.
This functional is analytically defined by the formula
$$
I=I(\Psi_{1},\dots, \Psi_{N})=\frac{\int\limits_{U_1}{\ldots \int\limits_{U_N}{A(u_{1},\ldots ,u_{N})d\Psi_{1}(u_{1})\ldots
d\Psi_{N}(u_{N})}}}{\int\limits_{U_1}{\ldots \int\limits_{U_N}{B(u_{1},\ldots ,u_{N})d\Psi_{1}(u_{1})\ldots
d\Psi_{N}(u_{N})}}},
\eqno(51)
$$
where the integrands in the numerator and denominator have the form
$$
A(u_{1},\ldots,u_{N})=\sum\limits_{i=1}^{N}{d_{i}(u_{i})}{\widehat{D}}^{(i)}
(u_{1},\ldots,u_{i-1},u_{i+1},\ldots,u_{N}),
\eqno(52)
$$
$$
B(u_{1},\ldots,u_{N})=\sum\limits_{i=1}^{N}{T_{i}(u_{i})}{\widehat{D}}^{(i)}
(u_{1},\ldots,u_{i-1},u_{i+1},\ldots,u_{N}).
\eqno(53)
$$
In turn, the functions ${\widehat{D}}^{(i)}(u_{1},\ldots,u_{i-1},u_{i+1}, \ldots ,u_{N})$,
$i=1,2,\dots,N$, contained in the right-hand sides of relations~(52) and~(53)
are defined by the formulas
$$
{\widehat{D}}^{(i)}(u_{1}, \ldots ,u_{i-1},u_{i+1}, \ldots , u_{N})=
$$
$$
=(-1)^{N+i+2}\sum\limits_{\alpha ^{(N),i}}{(-1)}^{\delta (\alpha
^{(N),i}) }{\widehat{D}}_{0}^{(i)}(\alpha ^{(N),i},u_{1}, \ldots
,u_{i-1},u_{i+1}, \ldots , u_{N}),
\eqno(54)
$$
where $\alpha ^{(N),i}=(\alpha_{1},\ldots,\alpha_{i-1},\alpha_{i+1},\ldots,\alpha _{N})$
is an arbitrary permutation of the numbers\break $(1,\ldots,i-1,i+1,\ldots,N)$,
$\delta(\alpha ^{(N),i})$ is the number of inversions in the permutation $\alpha^{(N),i}$, and
$$
{\widehat{D}}_{0}^{(i)}(\alpha^{(N),i},u_{1},\ldots,u_{i-1},u_{i+1},\ldots,u_{N})=
$$
$$
 ={\widetilde{p}}_{1,\alpha _{1}}(u_{1})\ldots {\widetilde{p}}_{i-1,\alpha
_{i-1}}(u_{i-1}){\widetilde{p}}_{i+1,\alpha _{i+1}}(u_{i+1})\ldots
{\widetilde{p}}_{N,\alpha _{N}}(u_{N}),
\eqno(55)
$$
$$
{\widetilde{p}}_{k,\alpha _{k}}(u_{k})=
\begin{cases}
p_{kk}(u_{k})-1,~ \alpha _{k}=k, \\
p_{k,\alpha _{k}}(u_{k}),~ \alpha _{k}\ne k, \\
\end{cases}
k=1, \ldots, i-1, i+1,\ldots,N.
\eqno(56)
$$
The functions $p_{ij}(u_i)$, $T_{i}(u_{i})$, and $d_{i}(u_{i})$, $u_i\in U_i$, $i,j=1,2,\dots,N$,
contained in (52)--(56) are defined by formulas~(42), (44), and~(48), respectively.}

Theorem~2 was proved in~\cite{23} in a very concise form.
The reader interested in a more detailed justification of this result
is referred to the text of A.V.~Ivanov's Candidate thesis (\cite{27}, Chapter~3).

So Theorem~2 permits obtaining an explicit analytic representation
of the stationary cost index of the form~(50)
as a liner-fractional integral functional of the set of probability measures
$\Psi=(\Psi_{1},\Psi_{2},\ldots,\Psi_{N})$ determining a control strategy for the semi-Markov process~$\xi(t)$.
In its analytic form, the functional $I(\Psi_1,\Psi_{2},\dots,\Psi_{N})$
given by~(51) completely coincides with the liner-fractional integral functional of the form~(4).
In this case, the integrands in the numerator and denominator are given by formulas~(52), (53)
and auxiliary relations (54)--(56). Thus, the test function
of the liner-fractional integral functional~(51),
$$
C(u_1, u_2,\dots, u_N)=\dfrac{A(u_1, u_2,\dots, u_N)}{B(u_1, u_2,\dots, u_N)},
\eqno(57)
$$
is also explicitly determined by formulas~(57), (52), and~(53).

We note that not only the stationary cost index of the form~(50)
depends as a liner-fractional integral
on the probability measures that form a control strategy.
Already in~\cite{24}, the author of this study proved the corresponding assertions
on the structure of several functionals related to the holding time of a semi-Markov process
in a given finite subset of states.
Thus, the statement of the optimal control problem for a semi-Markov process
in the form of the unconditional extremum problem
for a liner-fractional integral functional is completely meaningful.
\bigskip

\textbf{5. Theoretical solution of the optimal control problem.}

In this section we present the main result on the optimal control strategies for the semi-Markov process,
which is ideologically related to Theorem~1. To obtain a correct and meaningful statement of this result,
it is necessary to analyze several preliminary conditions under which the above-posed extremum problem
is considered. We perform such an analysis in the form of remarks concerned with several aspects of the problem
in question.
\medskip

\textbf{Remark~2.} Conditions~1--3 in the BPC system are necessary for the well-conditioned statement
of an unconditional extremum problem for a liner-fractional integral functional.
If this functional is the quality index in the optimal control problem for a stochastic process,
then it is necessary to supplement these conditions with an additional condition related to
some regularity of the decisiion process itself.
Namely, there must exist a specific meaningful index related to the behavior of this process
and be representable as a liner-fractional integral functional.
If it is required that the ergodic relation~(50) be satisfied,
then one must use Theorem~2 and formulate the optimal control problem
for the liner-fractional integral functional~(51) in the form~(5).
Thus, it is necessary to introduce a condition that ensures
the existence of a unique stationary distribution of the embedded Markov chain and
the satisfaction of relation~(50).
By analogy with~\cite{19}, Chapter~5, we formulate this additional condition as follows.

\begin{enumerate}
\item[4.] For any control strategy $\Psi=(\Psi_1,\Psi_2,\dots,\Psi_N)\in\Gamma$,
the embedded Markov chain of the semi-Markov process $\xi(t)$
has precisely one class of recursive positive states.
\end{enumerate}

In what follows, the BPC system is the set of conditions comprising conditions~1--3
formulated in Section~1 and additional condition~4.
\medskip

\textbf{Remark~3.}
Below we formulate and prove the main theorem on the optimal control strategy
for a semi-Markov process with finitely many states. We formulate this theorem with respect to
the extremum problem~(5), where the goal functional $I(\Psi_1, \Psi_2,\dots,\Psi_N)$
has the form of a general liner-fractional integral functional given by formula~(4).
This is related to the fact that the goal functional in the optimal control problem
need not have a character of stationary cost index of the form~(50).

As was already noted at the end of the preceding section,
there also exist other meaningful characteristics of control quality
whose dependence on the set of probability measures determining the control strategy
is liner-fractional.
Thus, the control problem posed above is more general than the problem
where the goal functional is a stationary cost index of the form~(50).

Now we introduce the notion of admissible control strategy
for a semi-Markov process with finitely many states.

\textbf{Definition~2.}
A control strategy $\Psi=(\Psi_1,\Psi_2,\dots,\Psi_N)$
is said to be admissible in this problem if it satisfies conditions~1--4
of the BPC system.
\medskip

\textbf{Remark~4.}
Remark~1 implies that if we demand that the function $B(u_1,u_2,\dots,u_N)$
be strictly of constant sign for all $(u_1, u_2,\dots,u_N)\in U$,
then we can assume that the strategies $(\Psi_1, \Psi_2,\dots,\Psi_N)$
satisfying conditions~1, 3, and~4 of the BPC system
are admissible.
With regard to Remark~7 below on the natural character of the condition
that the function $B(u_1,u_2,\dots,u_N)$ is strictly of constant sign for all values of the arguments
$(u_1, u_2,\dots,u_N)\in U$,
it is required that this condition is satisfied in the formulation
of the main theorem given below on the optimal control strategy for a semi-Markov process.
\medskip

\textbf{Remark~5.}
If we consider the optimal control problem for a semi-Markov process,
where the goal functional does not coincide with the stationary cost index~(50),
then we possibly need other additional conditions that ensure the existence
of this index and its representation in the form~(51).
Therefore, in the statement of the main theorem, we use the term
``admissible strategies'' in the wide sense under the assumption
that all necessary conditions are satisfied for each control quality index under study.
\medskip

\textbf{Remark~6.}
The set of admissible strategies need not coincide with the set of all possible
control strategies. In particular, the admissible strategies can consist
only of discrete probability measures $\Psi_1,\Psi_2,\dots,\Psi_N$,
i.e., of measures that are concentrated on arbitrary discrete sets of points
belonging to the spaces $U_1,U_2,\dots,U_N$.
\medskip

\textbf{Remark~7.}
If as the goal functional $I(\Psi_1,\Psi_2,\dots,\Psi_N)$ of the extremum problem~(5),
we consider the stationary cost index~(50) representable in the form~(51),
then the function $B(u_1,u_2,\dots,u_N)$ has the following theoretical meaning.
This function is the conditional mathe\-matical expectation
of the length of the time period between neighboring times of variation in the state
of the semi-Markov process $\xi(t)$ under the condition that its control strategy
is deterministic and is determined by a set of argument values $(u_1,u_2,\dots,u_N)$.
Then the condition that the function $B(u_1,u_2,\dots,u_N)$
is strictly of constant sign for all $(u_1,u_2,\dots,u_N)\in U$ is natural
and actually means that,
for any prescribed deterministic control strategy,
the process $\xi(t)$ has no instantaneous states for which the holding time is zero.
\medskip

\textbf{Remark 8.}
We make several notes related to the integrand in the numerator
of liner-fractional integral functional~(51).
As previously, as the goal functional $I(\Psi_1, \Psi_2,\dots,\Psi_N)$
of the extremum problem~(5), we consider the stationary cost index~(50).
Then for any fixed set of argument values $(u_1,u_2,\dots,u_N)\in U$,
the value of the function $A(u_1,u_2,\dots,u_N)$
is the conditional mathematical expectation of the increment of the cost functional
during the time of stay of the semi-Markov process $\xi(t)$
in a fixed state provided that the control strategy is deterministic
and is determined by the set $(u_1,u_2,\dots,u_N)\in U$ mentioned above.
We note that, in the theorem on the extremum of a liner-fractional integral functional,
which was proved in~\cite{6} (Chapter~10),
the integrand in the numerator is assumed to be bounded on the whole set of argument values.
For many mathematical models and the optimal control problems related to them,
this condition is too restrictive.
As examples, we can consider models of optimal continuous-product inventory control
studied in~\cite{28},~\cite{29}.
In this study, we impose the only condition that the function $A(u_1,u_2,\dots,u_N)$
is integrable with respect to any given set of probability measures
$\Psi=(\Psi_1, \Psi_2,\dots,\Psi_N)$ which forms a control strategy for the semi-Markov process~$\xi(t)$
(condition~1 in the BPC system).

Now we formulate and prove the main theorem on the existence and the form of representation of
the optimal control strategy for a semi-Markov process with finitely many states.

\textbf{Theorem 3.}
{\it Consider the optimal control problem for the semi-Markov process
$\xi(t)$ in the form of the extremum problem~(5) which is defined on the set of admissible strategies~$\Gamma$
for a liner-fractional integral functional of general form~(4).
Let the function $B(u_1,u_2,\dots,u_N)$ in the definition of the functional~(4)
be strictly of constant sign (strictly positive or strictly negative)
for all argument values $(u_1,u_2,\dots,u_N)\in U$.
Then the following assertions hold.
\begin{enumerate}
\item[{\rm(1)}]
If the test function $C(u_1,u_2,\dots,u_N)=\dfrac{A(u_1,u_2,\dots,u_N)}{B(u_1,u_2,\dots,u_N)}$
is bounded above or below and attains a global extremum on the set
$U=U_1\times U_2\times\dots\times U_N$ (maximum or minimum),
then the optimal control strategy for a semi-Markov process exists and is attained
on the deterministic control strategy
$\Psi^{*(+)}=\left(\Psi_1^{*(+)},\Psi_2^{*(+)},\dots,\Psi_N^{*(+)}\right)\in\Gamma^*$
or on the deterministic control strategy
$\Psi^{*(-)}=\left(\Psi_1^{*(-)},\Psi_2^{*(-)},\dots,\Psi_N^{*(-)}\right)$
depending on the form of the extremum (maximum or minimum).
In this case, the degenerate probability measures\break
$\Psi_1^{*(+)},\Psi_2^{*(+)},\dots,\Psi_N^{*(+)}$ are concentrated at the respective points
$u_1^{(+)},u_2^{(+)},\dots,u_N^{(+)}$
if $u^{(+)}=\left(u_1^{(+)},u_2^{(+)},\dots,u_N^{(+)}\right)\in U$ is the point of global maximum
of the function $C(u_1,u_2,\dots,u_N)$,
and then the following relations are satisfied:
$$
\max\limits_{\Psi \in \Gamma} I(\Psi)=\max\limits_{\Psi_i \in \Gamma_i, i=\overline{1,N}} I(\Psi_1,\Psi_2,\dots,\Psi_N)=
\max\limits_{\Psi_i^* \in \Gamma_i^*, i=\overline{1,N}} I(\Psi_1^*,\Psi_2^*,\dots,\Psi_N^*)=
$$
$$
=I\left(\Psi_1^{*(+)},\Psi_2^{*(+)},\dots,\Psi_N^{*(+)}\right)
=\max\limits_{(u_1,u_2,\dots,u_N)\in U}\dfrac{A(u_1,u_2,\dots,u_N)}{B(u_1,u_2,\dots,u_N)}=
$$
$$
=\dfrac{A\left(u_1^{(+)}, u_2^{(+)},\dots,u_N^{(+)}\right)}{B\left(u_1^{(+)}, u_2^{(+)},\dots,u_N^{(+)}\right)}.
\eqno(58)
$$
Similarly, the degenerate probability measures $\Psi_1^{*(-)},\Psi_2^{*(-)},\dots,\Psi_N^{*(-)}$
are concentrated at the respective points $\left(u_1^{(-)}, u_2^{(-)},\dots,u_N^{(-)}\right)$
if $u^{(-)}=\left(u_1^{(-)},u_2^{(-)},\dots,u_N^{(-)}\right)\in U$ is the point of global minimum
of the function $C(u_1,u_2,\dots,u_N)$,
and then the following relations are satisfied:
$$
\min\limits_{\Psi \in \Gamma} I(\Psi)=\min\limits_{\Psi_i \in \Gamma_i, i=\overline{1,N}} I(\Psi_1,\Psi_2,\dots,\Psi_N)=
\min\limits_{\Psi_i^* \in \Gamma_i^*, i=\overline{1,N}} I(\Psi_1^*,\Psi_2^*,\dots,\Psi_N^*)=
$$
$$
=I\left(\Psi_1^{*(-)},\Psi_2^{*(-)},\dots,\Psi_N^{*(-)}\right)
=\min\limits_{(u_1,u_2,\dots,u_N)\in U}\dfrac{A(u_1,u_2,\dots,u_N)}{B(u_1,u_2,\dots,u_N)}=
$$
$$
=\dfrac{A\left(u_1^{(-)}, u_2^{(-)},\dots,u_N^{(-)}\right)}{B\left(u_1^{(-)}, u_2^{(-)},\dots,u_N^{(-)}\right)}.
\eqno(59)
$$
\item[{\rm(2)}]
If the test function $C(u_1,u_2,\dots,u_N)=\dfrac{A(u_1,u_2,\dots,u_N)}{B(u_1,u_2,\dots,u_N)}$
is bounded above or below but does not attain the global extremum on the set
$U=U_1\times U_2\times\dots\times U_N$,
then, for any $\varepsilon>0$, one can choose
an $\varepsilon$-optimal deterministic control strategy
$\Psi^{*(+)}(\varepsilon)=$\break $=\left(\Psi_1^{*(+)}(\varepsilon),\Psi_2^{*(+)}(\varepsilon),\dots,\Psi_N^{*(+)}(\varepsilon)\right)$
or a similar $\varepsilon$-optimal deterministic control strategy
$\Psi^{*(-)}(\varepsilon)=\left(\Psi_1^{*(-)}(\varepsilon),\Psi_2^{*(-)}(\varepsilon),\dots,\Psi_N^{*(-)}(\varepsilon)\right)$
depending on the form of the extremum (maximum or minimum).
In this case,
the degenerate probability measures\break
$\Psi_1^{*(+)}(\varepsilon),\Psi_2^{*(+)}(\varepsilon),\dots,\Psi_N^{*(+)}(\varepsilon)$
are concentrated at the respective points\break
$u_1^{(+)}(\varepsilon),u_2^{(+)}(\varepsilon),\dots,u_N^{(+)}(\varepsilon)$, where
$u^{(+)}(\varepsilon)=\left(u_1^{(+)}(\varepsilon),u_2^{(+)}(\varepsilon),\dots,u_N^{(+)}(\varepsilon)\right)\in U$
is an arbitrary point such that
$$
\sup\limits_{(u_1,u_2,\dots,u_N) \in U}\dfrac{A(u_1,u_2,\dots,u_N)}{B(u_1,u_2,\dots,u_N)}-\varepsilon <
\dfrac{A\left(u_1^{(+)}(\varepsilon),u_2^{(+)}(\varepsilon),\dots,u_N^{(+)}(\varepsilon)\right)}
{B\left(u_1^{(+)}(\varepsilon),u_2^{(+)}(\varepsilon),\dots,u_N^{(+)}(\varepsilon)\right)}<
$$
$$
<\sup\limits_{(u_1,u_2,\dots,u_N) \in U}\dfrac{A(u_1,u_2,\dots,u_N)}{B(u_1,u_2,\dots,u_N)}<\infty
\eqno(60)
$$
if the function $C(u_1,u_2,\dots,u_N)$ is bounded above and extremum problem~(5)
is the maximum problem.
Similarly, the probability measures
$\Psi_1^{*(-)}(\varepsilon),\Psi_2^{*(-)}(\varepsilon),\dots,\Psi_N^{*(-)}(\varepsilon)$
are concentrated at the respective points
$u_1^{(-)}(\varepsilon),u_2^{(-)}(\varepsilon),\dots,u_N^{(-)}(\varepsilon)$, where
$u^{(-)}(\varepsilon)=\left(u_1^{(-)}(\varepsilon),u_2^{(-)}(\varepsilon),\dots\right.$,
\break $\left.u_N^{(-)}(\varepsilon)\right)\in U$
is an arbitrary point such that
$$
-\infty<\inf\limits_{(u_1,u_2,\dots,u_N) \in U}\dfrac{A(u_1,u_2,\dots,u_N)}{B(u_1,u_2,\dots,u_N)} <
\dfrac{A\left(u_1^{(-)}(\varepsilon),u_2^{(-)}(\varepsilon),\dots,u_N^{(-)}(\varepsilon)\right)}
{B\left(u_1^{(-)}(\varepsilon),u_2^{(-)}(\varepsilon),\dots,u_N^{(-)}(\varepsilon)\right)}<
$$
$$
<\inf\limits_{(u_1,u_2,\dots,u_N) \in U}\dfrac{A(u_1,u_2,\dots,u_N)}{B(u_1,u_2,\dots,u_N)}+\varepsilon
\eqno(61)
$$
if the function $C(u_1,u_2,\dots,u_N)$ is bounded below and the extremum problem~(5)
is the minimum problem.
\item[\rm{(3)}]
If the test function $C(u_1,u_2,\dots,u_N)=\dfrac{A(u_1,u_2,\dots,u_N)}{B(u_1,u_2,\dots,u_N)}$
is not bounded above or below,
then there exists no optimal control strategy in the sense of the corresponding extremum problem.
In this case, there exists a sequence of deterministic control strategies
$\Psi^*(n)=\left(\Psi_1^{*(+)}(n),\Psi_2^{*(+)}(n),\dots,\Psi_N^{*(+)}(n)\right)$, $n=1,2,\dots$,
where the degenerate probability measures $\Psi_1^{*(+)}(n),\Psi_2^{*(+)}(n),\dots,\Psi_N^{*(+)}(n)$
are concentrated at the respective points\break
$u_1^{(+)}(n),u_2^{(+)}(n),\dots,u_N^{(+)}(n)$, $n=1,2,\dots$,
for which the following relations are satisfied:
$$
I\left(\Psi^{*(+)}(n)\right)=I\left(\Psi_1^{*(+)}(n),\Psi_2^{*(+)}(n),\dots,\Psi_N^{*(+)}(n)\right)=
$$
$$
=\dfrac{A\left(u_1^{(+)}(n),u_2^{(+)}(n),\dots,u_N^{(+)}(n)\right)}{B\left(u_1^{(+)}(n),u_2^{(+)}(n),\dots,u_N^{(+)}(n)\right)}
\to \infty~\mbox{ as } n\to\infty
\eqno(62)
$$
if the function $C(u_1,u_2,\dots,u_N)$ is not bounded above.
Similarly,
there exists a sequence of deterministic control strategies
$\Psi^{*(-)}(n)=\left(\Psi_1^{*(-)}(n),\Psi_2^{*(-)}(n),\dots,\Psi_N^{*(-)}(n)\right)$, $n=1,2,\dots$,
where the degenerate probability measures $\Psi_1^{*(-)}(n),\Psi_2^{*(-)}(n),\dots,\Psi_N^{*(-)}(n)$
are concentrated at the respective points
$u_1^{(-)}(n),u_2^{(-)}(n),\dots,u_N^{(-)}(n)$, $n=1,2,\dots$,
for which the following relations are satisfied:
$$
I\left(\Psi^{*(-)}(n)\right)=I\left(\Psi_1^{*(-)}(n),\Psi_2^{*(-)}(n),\dots,\Psi_N^{*(-)}(n)\right)=
$$
$$
=\dfrac{A\left(u_1^{(-)}(n),u_2^{(-)}(n),\dots,u_N^{(-)}(n)\right)}{B\left(u_1^{(-)}(n),u_2^{(-)}(n),\dots,u_N^{(-)}(n)\right)}
\to -\infty~\mbox{ as } n\to\infty
\eqno(63)
$$
if the function $C(u_1,u_2,\dots,u_N)$ is not bounded below.
\end{enumerate}
In this case, the above-formulated assertions in each item of Theorem~3
can be satisfied either separately for one of the two types of extremum
or jointly for both types of extremum.}
\medskip

\textbf{Proof of Theorem~3.} If the conditions of Theorem~3 are satisfied,
then the conditions of Theorem~1 are also satisfied.
Thus, to solve the optimal control problem for a semi-Markov process
in the form of the extremum problem~(5) with a liner-fractional integral functional of the form~(4),
one can use Theorem~1.
Then all assertions of Theorem~3 directly follow from the corresponding assertions of Theorem~1.
In particular, analytic relations (58)--(63) are consequences of relations (6)--(11)
contained in the statement of Theorem~1.
\bigskip

\textbf{6. Application of obtained theoretical results.}

Theorems~2 and~3 proved above form the theoretical foundation of the new method for studying
optimal control problems for a semi-Markov process with a finite set of states.
In the whole, this method consists
in determining an explicit analytic representation for the test function
of the liner-fractional integral functional which plays the role of the control quality index
and in using the assertions of Theorem~3 that determine the existence and
the character of the optimal control strategy.
Note that such a study of problems of optimal control of stochastic systems
was in fact already performed in several works of P.V.~Shnurkov and his coauthors.
In particular, in~\cite{30}, a model of control was considered for
a terminating renewal process which describes the operation and maintenance of a technical system
on the period of its operation till the first failure.
The control problem was solved for different values
of efficiency and reliability indices of this system
which have the structure of a liner-fractional integral functional.
We note that qualitatively different results related to the character of optimal control,
which correspond to different assertions of Theorem~3,
were already obtained in this paper.

Models of regenerative processes used to study inventory control systems
were considered in~\cite{28},~\cite{29}.
Different control quality indices were represented as liner-fractional integral functionals.
The test functions of these functionals were obtained explicitly
and were investigated for the global extremum.
A sufficiently complicated semi-Markov process model with finitely many states,
which describes the continuous-product inventory control system,
was considered in~\cite{23}, \cite{31}.
In this model, the control quality indices had also had the structure
of liner-fractional integral functionals, and explicit analytic representations
were obtained for their test functions.
We also mention~\cite{32},~\cite{33}, where the semi-Markov model with discretely-continuous phase space
was studied. The control quality indices in this model were obtained explicitly
as functions of two continuous parameters of control.

In all above-listed works, the method used
to solve the optimal control problem for regenerative or semi-Markov stochastic processes
is in fact based on an analysis of the extremum properties of the test function
of the corresponding liner-fractional integral functional.
The arguments used in Section~2 to analyze the results of previous studies
on the theory of control of a semi-Markov processes with finitely many states
show that this method had not been rigorously justified when these works were written and published.
But after the publication of~\cite{10} and the present study,
one can assert that the results obtained there are completely theoretically justified.
\bigskip

\textbf{Conclusion.}

Let us summarize the above investigation.
For this, we briefly describe some provisions characterizing the novelty
and theoretical value of the obtained results.

{\bf1.} In this paper, sufficient conditions for the existence of
a deterministic optimal control strategy
for a semi-Markov process with finitely many states are obtained.
In recent studies~\cite{21}, \cite{22},
it was proved that, under very general assumptions on the state and decision spaces,
the optimal control strategy for a semi-Markov process is deterministic and stationary.
The novelty and the main specific characteristic of the results obtained in this paper
is that the extremum problem for the liner-fractional integral functional
and the optimal control problem for a semi-Markov stochastic process
can be solved by studying the extremum properties of the test function
$C(u_1,u_2,\dots,u_N)=\dfrac{A(u_1,u_2,\dots,u_N)}{B(u_1,u_2,\dots,u_N)}$.
If the test function attains the global extremum at a certain point,
then the solution of the corresponding optimal control problem for the semi-Markov process
exists and is attained on the deterministic control strategy
determined by the point of global extremum.
Thus, determining the global extremum of the test function
permits simultaneously solving the existence problem of optimal control strategy
and determining the optimal strategy itself.
We note that, in many specific problems, the test function can be expressed explicitly
and its properties can be studies analytically of numerically
(see, e.g.,~\cite{30},~\cite{28}, \cite{29},~\cite{31},~\cite{23} cited above).
The author of the present study also proved general assertions
about the analytic representations of test functions of different types of meaningful indices
related to the control of semi-Markov processes
(\cite{23},~\cite{24}).

{\bf2.} The assertions obtained above also give an exhaustive answer
to the question of the solution of the extremum problem
for liner-fractional integral functional and
of the optimal control problem for a semi-Markov process
in the case where the test function
$C(u_1,u_2,\dots,u_N)=\dfrac{A(u_1,u_2,\dots,u_N)}{B(u_1,u_2,\dots,u_N)}$
does not attain the the global extremum.
The following two situations are possible in this case.
If the function $C(u_1,u_2,\dots,u_N)$ is bounded above or below
on the set $U=U_1\times U_2 \times \dots \times U_N$,
i.e., if the set of its values has a finite upper or lower bound,
then there exists a degenerate probability measure on which
the goal functional takes a value lying in an arbitrarily small neighborhood
of its upper or lower bound.
As for the optimal control problem for a semi-Markov process,
this assertion can be interpreted in the sense that there exists
a deterministic $\varepsilon$-optimal control strategy.
But if the function $C(u_1,u_2,\dots,u_N)$ is not bounded above or below,
then the set of solutions of the liner-fractional integral goal functional
is also unbounded.
There does not exist a solution of the extremum problem for such a functional
(the maximum or minimum problem).
Then, in the framework of the above-posed problem,
there also does not exist an optimal control strategy
for a semi-Markov stochastic process.

Thus, the above-proved theorem on the unconditional extremum
of a liner-fractional integral functional gives an exhausting solution
for the whole class of optimal control problems for semi-Markov stochastic processes
with finitely many states.
\bigskip

\textbf{Acknowledgments.}

The author wishes to express deep gratitude
to the eminent scientist,
academician of the National Academy of Sciences of Ukraine,
Professor Vladimir Semenovich Korolyuk
for useful discussions of the obtained results
and invaluable personal moral support during many years.


\begin{thebibliography}{99}
\bibitem{10}
{\em Shnurkov~P.\,V.},
``Solution of the unconditional extremum problem
for a linear-fractional integral functional on a set of probability measures'',
Dokl. Ross. Akad. Nauk \textbf{94} (2), 550--554 (2016).
[Dokl. Math. \textbf{470} (4), 387--392].

\bibitem{1}
{\em Bajalinov E.\,B.},
{\em Linear-Fractional Programming. Theory, Methods, Applications and Software}
(Kluwer, Boston--Dordrecht--London, 2003).

\bibitem{2}
{\em Borza M., Rambely A.S., and Saraj M.},
``Solving linear fractional programming problems with interval coefficients
in the objective function. A new approach'',
Appl. Math. Sci. Ruse \textbf{6} (69), 3443--3452 (2012).

\bibitem{3}
{\em Joshi V.D., Singh E., and Gupta N.},
``Primal-dual approach to solve linear fractional programming problem'',
J. Appl. Math., Statist. Inform. (JAMSI) \textbf{4} (1), 61--69 (2008).

\bibitem{4}
{\em Hasan M.B. and Acharjee S.},
``Solving LFP by converting it into a single LP,''
Int. J. Oper. Res. \textbf{8} (1), 1--14 (2011).

\bibitem{5}
{\em Barzilovich E.\,Yu. and Kashtanov V.\,A.},
{\em Several Mathematical Problems of Queuing Theory for Complex Systems}
(Sov. Radio, Moscow, 1971)
[in Russian].

\bibitem{6}
{\em Problems of Mathematical Reliability Theory},
ed. by B.\,V.~Gnedenko
(Radio i Svyaz', Moscow, 1983)
[in Russian].

\bibitem{7}
{\em Shiryaev A.\,N.},
{\em Probability}
(MTzNMO, Moscow, 2011)
[in Russian].

\bibitem{8}
{\em Hennequin~P.\,L. and Tortrat~A.},
{\em Th\'eorie des probabilit\'es et quelques applications}
(Masson, Paris, 1965; Nauka, Moscow, 1974).

\bibitem{9}
{\em Halmos P.}
{\em Measure Theory}
(Inostrannaya Literatura, Moscow, 1953),
[Russian translation].

\bibitem{11}
{\em Kudryavtsev L.\,D.},
{\em Course of Mathematical Analysis}
(Drofa, Moscow, 2006), Vol.~1
[in Russian].

\bibitem{12}
{\em Howard, R.\,A.},
{\em Dynamic Programming and Markov Processes}
(MIT Press, Cambridge, Massachusetts, 1960; Sov. Radio, Moscow, 1964).

\bibitem{13}
{\em Rykov~V.\,V.},
``Markov decision processes with finite state and decision spaces'',
Teor. Veroyatnost. Primenen. \textbf{11} (2), 343--351 (1966)
[Theory Probab. Appl. \textbf{11}, 302--311 (1966)].

\bibitem{14}
{\em Jewell W.\,S.},
``Semi-Markov decision processes'',
in {\em Cybernetic Collection}
(Mir, Moscow, 1967), No.~4, pp.~97--134.

\bibitem{15}
{\em Fox B.},
``Markov renewal programming by linear fractional programming'',
SIAM J. Appl. Math. \textbf{14}, 1418--1432 (1996).

\bibitem{16}
{\em Denardo E.\,V.},
``Contraction mappings in the theory underlying dynamic programming'',
SIAM Rev. \textbf{9}, 165--177 (1967).

\bibitem{17}
{\em Howard R.\,A.},
``Research in semi-Markovian decision structures'',
J. Oper. Res. Soc. Japan \textbf{6}, 163--199 (1963).

\bibitem{18}
{\em Osaki S. and Mine H.},
``Linear programming algorithms for Markovian decision processes'',
J. Math. Anal. Appl. \textbf{22}, 356--381 (1968).

\bibitem{19}
{\em Mine H. and Osaki S.},
{\em Markovian Decision Processes}
(Elsevier, Amsterdam, 1970; Nauka, Moscow, 1977).

\bibitem{20}
{\em Gikhman I.\,I. and Skorokhod A.\,V.},
{\em Controlled Stochastic Processes}
(Naukova Dumka, Kiev, 1977)
[in Russian].

\bibitem{21}
{\em Luque-Vasquez F. and Hernandez-Lerma O.},
``Semi-Markov control models with average costs'',
Appl. Math. \textbf{26} (3), 315--331 (1999).

\bibitem{22}
{\em Vega-Amaya O. and Luque-Vasquez F.},
``Sample-path average cost optimality for semi-Markov control processes
on Borel spaces: unbounded costs and mean holding times'',
Appl. Math. \textbf{27} (3), 343--367 (2000).

\bibitem{25}
{\em Korolyuk V.\,S. and Turbin A.\,F.},
{\em Semi-Markov Process and Their Applications}
(Naukova Dumka, Kiev, 1976)
[in Russian].

\bibitem{26}
{\em Janssen J. and Manca R.},
{\em Applied Semi-Markov Processes}
(Springer, New York, 2006)

\bibitem{23}
{\em Shnurkov~P.\,V. and Ivanov~A.\,V.},
``Analysis of a discrete semi-Markov model of continuous inventory control
with periodic interruptions of consumption'',
Diskret. Mat. \textbf{26} (1), 143--154 (2014)
[Discrete Math. Appl. \textbf{25} (1), 59--67 (2015)].

\bibitem{27}
{\em Ivanov~A.\,V.},
{\em Analysis of Discrete Semi-Markov Model of Continuous Inventory Control
with Periodic Interruptions of Consumption},
Cand. Sci. (Phys.--Math.) Dissertation
(NRU ``Higher School of Economics'', Moscow, 2014)
[in Russian].

\bibitem{24}
{\em Shnurkov~P.\,V.},
{\em Methods for Studying Problems of Optimal Maintenance in Mathematical Reliability Theory},
Cand. Sci. (Phys.--Math.) Dissertation
(Moscow Institute of Electronic Engineering (MIEM), Moscow, 1983).

\bibitem{28}
{\em Shnurkov~P.\,V. and Mel'nikov R.\,V.},
``Optimal continuous-product inventory control in a regeneration model'',
in {\em Review of Applied and Industrial Mathematics}
\textbf{13} (3), 434--452 (2006).

\bibitem{29}
{\em Shnurkov~P.\,V. and Mel'nikov R.\,V.}
``Analysis of the problem of continuous-product inventory control under deterministic lead time'',
Avtomat. Telemekh. \textbf{10}, 93--113 (2008)
[Automat. Remote Control \textbf{69} (10), 1734--1751 (2008)].

\bibitem{30}
{\em Shnurkov~P.\,V.}
``Optimal maintenance of a technical system on the period till the first failure'',
in {\em Application of Analytic Methods in Probability Problems},
Collection of Scientific Papers
(Institute of Mathematics, Acad. Sci. UkrSSR, Kiev, 1986), pp~121--129.

\bibitem{31} {\em Shnurkov~P.\,V. and Ivanov~A.\,V.},
``Study of optimization problem in discrete semi-Markov model
of continuous-product inventory control'',
Vestnik Bauman Moscow State Technical Univ., Ser. ``Natural Sciences''
\textbf{3} (50), 63--87 (2013).

\bibitem{32}
{\em Shnourkoff P.\,V.},
``The two-element system with one restoring device optimum maintenance'',
Stochastic Anal. Appl. \textbf{15} (5), 823--837 (1997).

\bibitem{33}
{\em Shnourkoff P.\,V.},
``The two-element system optimum maintenance tills the first fail'',
Stochastic Anal. Appl. \textbf{19} (6), 1005--1024 (2001).

\end{thebibliography}
\end{document}